\documentclass[11pt]{article}
\usepackage{latexsym} 
\usepackage{amsfonts}
\usepackage{amsmath}
\usepackage{amsthm}
\usepackage{amssymb}
\renewcommand{\theequation}{\thesection.\arabic{equation}}
\newcounter{subequation}[equation]
\makeatletter

\newcommand{\p}{^{\prime}}

\expandafter\let\expandafter\reset@font\csname reset@font\endcsname

\def\subeqnarray{\arraycolsep1pt
    \def\@eqnnum\stepcounter##1{\stepcounter{subequation}%
        {\reset@font\rm(\theequation\alph{subequation})}}
\jot5mm     \eqnarray}

\makeatother

\def\be{\begin{equation}}
\def\ee{\end{equation}}
\def\bea{\begin{eqnarray}}
\def\eea{\end{eqnarray}}
\def\ba{\begin{array}}
\def\ea{\end{array}}
\def\dd{\partial}

\def\half{\frac{1}{2}}

\def\one#1{#1^{\raise5pt\hbox{$\scriptstyle\!\!\!\!1$}}\,{}}
\def\two#1{#1^{\raise5pt\hbox{$\scriptstyle\!\!\!\!2$}}\,{}}

\def\tilde{\widetilde}
\def\II{\hbox{{1}\kern-.25em\hbox{l}}}
\makeatletter
\def\binrel@#1{\begingroup
  \setboxz@h{\thinmuskip0mu
    \medmuskip\m@ne mu\thickmuskip\@ne mu
    \setbox\tw@\hbox{$#1\m@th$}\kern-\wd\tw@
    ${}#1{}\m@th$}%
  \edef\@tempa{\endgroup\let\noexpand\binrel@@
    \ifdim\wdz@<\z@ \mathbin
    \else\ifdim\wdz@>\z@ \mathrel
    \else \relax\fi\fi}%
  \@tempa
}
\let\binrel@@\relax
\def\overset#1#2{\binrel@{#2}%
  \binrel@@{\mathop{\kern\z@#2}\limits^{#1}}}
\def\underset#1#2{\binrel@{#2}%
  \binrel@@{\mathop{\kern\z@#2}\limits_{#1}}}
\makeatother
\newfont{\bbd}{msbm10 scaled\magstep1}

\newtheorem{prop}{Proposition}

\begin{document}

\begin{center}
{\LARGE {Intertwiners of Yangian representations 
}} 

\vspace{0.5cm}

{\large \sf

R. Kirschner\footnote{\sc e-mail:Roland.Kirschner@itp.uni-leipzig.de}} 

\vspace{0.5cm}

{\it Institut f\"ur Theoretische
Physik, Universit\"at Leipzig, \\
PF 100 920, D-04009 Leipzig, Germany}

\end{center}

\vspace{.5cm}
\begin{abstract}
\noindent
We construct type $g\ell(n)$ Yangian algebra evaluations of order $N$
embedded in Heisenberg algebras and consider their representations having a
highest weight. These Yangian algebra presentations depend on $nN$
parameters. We construct explicitly intertwiners, which exist if the
parameter arrays are related by permutations. 
The intertwiners are products of elementary adjacent parameter permutation
operators derived from $R$ operators obeying Yang-Baxter relations.  
Permutation coefficients  appear in the action on
representations. Their dependence on the parameters allows to distinguish
the types of representations. 

\end{abstract}

\renewcommand{\refname}{References.}
\renewcommand{\thefootnote}{\arabic{footnote}}
\setcounter{footnote}{0}

\section{Introduction}

It is well known that generators obeying the $g\ell(n)$ algebra relations can
be built of Heisenberg canonical pairs \cite{PJ35, JS, HP40}. 
This Jordan-Schwinger (JS) presentation allows to describe a 
one-parameter family
of  $g\ell(n)$ algebra representations. 
The $g\ell(n)$ algebra relations can be written in the form the Yang-Baxter
relation with Yang's R matrix and L operators composed of the 
latter generators. The explicit form of the L operators allows to construct
R operators intertwining JS representations. A particular form of these
R operators in terms of a contour integral over a shift operator
appeared in the context of scattering amplitude calculations
\cite{DHP09,AH10,CDK13,FK17,RK23}.

The JS presentation of the  $g\ell(n)$ algebra  does not cover generic
representations. The Biedenharn iterative construction \cite{BL92,BL95} 
allows to 
obtain the generators  for generic representations starting from
the JS form of generators. We explain the L operator formulation of
this construction \cite{KK09,RK13}. It starts from products of JS type L
operators. Constraints are imposed on the involved canonical pairs to obtain 
the first order evaluation form. We rely on Dirac's Hamilton formalism
\cite{PD}.
In the framework of this formulation we construct operators of
representation parameter
permutation from R operators by imposing the latter constraints.
The operators of parameter permutation  obtained in this way are compatible
with the ones obtained  on the basis of the method of factorization
\cite{SD05,DKK07,DKKV08}, which appears more involved at higher rank.  

We consider the order N evaluations of the Yangian algebra in terms of the
monodomies of the N-fold product of JS type L operators and  their
representations with a highest weight. 
The  R operators intertwining L operators of JS type allow to obtain the
intertwining operators of
Yangian algebra evaluations with parameters related by permutations.
The parameter permutation
operators intertwine different presentations with the coinciding center. 
Intertwining operators exist only if their
representation parameters are related by permutations. 

We observe features generally well known in the 
representation theory of Lie algebras
\cite{BGG71}. We formulate the corresponding relations for the 
order N evaluation of the $g\ell(n)$ type Yangian algebra.   
Yang-Baxter relations are a main working tool. Starting with the
Jordan-Schwinger form we are able to provide explicit expressions
of the relevant Yang-Baxter and the  intertwining operators. 

The representations are expressed in terms of homogeneous functions of the
coordinates out of the underlying Heisenberg canonical pairs. The parameters
appear as degrees of homogeneity. Particular functions represent the highest
weight vectors and the representations are spanned by monomials multiplying
them. The intertwining operators map the highest weight functions in
particular. Quite involved functions may appear in this way out of simple
ones. In a particular case of one-dimensional representations  highest
weight functions, called Yangian symmetric correlators, 
are related to scattering amplitudes of gauge field theories and
to kernels describing the high-energy asymptotics of scattering
\cite{CDK13,RK23}.

In the action of the intertwining operators on the 
highest weight vectors  permutation coefficients
appear. They depend on the parameters taking complex  values in general.
Zeros and poles appear in the case of non-equivalent representations.
The permutation coefficients provide an alternative 
 way to
analyse the Yangian representations having a highest weight
and their reducibility.

We present the rank one case in details and show in particular how the  results
about Yangian representations composed of finite-dimensional $g\ell(2)$
representations \cite{NT97}, \cite{AM02,AM07}
can be reproduced in analyzing the permutation coefficients.
We construct
particular intertwining operators between  order four
evaluations of the Yangian algebra composed in equivalent way of two order
two evaluations. The parameter dependence of the related
permutation coefficient allows to distinguish the types of representations 
of the order two evaluations.

We consider beta sequences of permutations where the action of 
the involved permutation operators is easily defined. 
The highest weight functions appearing here are expressed in terms of
parameter dependent powers of
determinants composed of the normal coordinate vectors. The elementary
adjacent permutation operators act in a simple way resulting in determinant
factors of this type multiplied by Euler beta functions of the parameters.

We give explicit expressions of beta sequence operators in the order two
evaluation for the parameter pair permutations of any parameter of a factor
 with the adjacent parameter of the adjacent
factor. 
 Using these we
construct the beta sequence  for the parameter pair permutations
of the parameters at position $i$ of the first and the second factors
 distinguished by the coincidence of
the highest weight vectors. Therefore the corresponding permutation
coefficients play an important role in characterizing the type of
representations.

The paper is organized as follows:

The Yang-Baxter formulation of $g\ell(n)$ type Yangians is recalled in sect.
2 and the Jordan-Schwinger presentation is recalled in sect. 3. We consider
$R$ operators obeying Yang-Baxter relations with the Jordan-Schwinger type $L$
matrices in sect. 4. The Yang-Baxter formulation of the Biedenharn
construction is  presented in sect 5.  The monodromies composed of $N$
generic $L$ matrices are studied in sect. 6 with focus on the parameter
dependence.  In sect. 7 we consider the morphisms between 
Yangian presentations related by parameter permutations and obtain functions
representing highest weight vectors. In the $g\ell(2)$ case explicit expressions and details 
about 
the characterization of representations types by the parameter permutation
coefficients are given in sect. 8. In the general $g\ell(n)$ case 
the construction of beta sequences of 
permutation operators and the related permutation coefficients are given in
sect. 9.

\section{ Yangians of type $g\ell(n)$ }
\setcounter{equation}{0}

We recall well known notions and relations and refer to \cite{AM07} for more details.

We start with the formulation of the $g\ell(n)$ algebra in terms of the
fundamental Yang-Baxter relation involving Yang's R matrix,
$$ \mathcal{R}^{a_1a_2}_{\ \ b_1b_2} (u) = u \delta^{a_1}_{b_1} \delta^{a_2}_{b_2} + 
\delta^{a_1}_{b_2} \delta^{a_2}_{b_1} $$
and the L operator matrix with the $g\ell(n)$ algebra generators as matrix
elements,
$$ L^a_{\ \ b}(u) = u \delta^{a}_{b} + L^a_{\ \ b}. $$
\be \label{fundYB}
 \mathcal{R}(u-v) (L(u)\otimes I) (I\otimes L(v)) = (I\otimes L(v))( L(u)\otimes
\mathcal{R}(u-v),
\ee
$$ \mathcal{R}^{a_1a_2}_{\ \ b_1 b_2} (u-v) L^{b_1}_{\ \ c_1}(u) L^{b_2}_{\
\ c_2}(v) = L^{a_1}_{\ \ b_1}(v) L^{a_2}_{\ \ b_2}(u) \mathcal{R}^{b_1 b_2}_{\ \ c_1 c_2}
(u-v).
$$
The extended Yangian algebra can be defined by generators
\be \label{Tu} T^a_{\ \ b} (u) = 1 + \sum T^{[k]} u^{-k} \ee
obeying the analogon of (\ref{fundYB}) with  $ L^a_{\ \ b}(u)$ replaced by $T^a_{\ \ b}
(u)$. 

Let us consider Yangian  evaluations   of finite order $N$
where the monodromy $T(u)$ is constructed from the $L$ matrices in the way well know in 
the context of integrable models. 
Let the matrix elements of $L^I(u^I)$ act on the representation space $V^I$
and consider  the matrix product defining the monodromy matrix
\be \label{monodromy}
T(u,\delta^1, \dots,\delta^N)=  L^1(u^1) \cdots L^N(u^N), \quad
u^I = u+\delta^I,  \ee 
acting on the tensor product $ V^1 \otimes \cdots \otimes V^N $. The Yangian
co-multiplication property appears in the  well-known form:
If (\ref{fundYB}) holds then the monodromy $T(u)$ obeys the analogy of (\ref{fundYB}) 
with $L(u)$ substituted by $T(u)$. Comparing to (\ref{Tu}) here we 
restrict to representations where the higher generators $T^{[k]}, k>N $ act
as multiplication by zero and redefine by multiplying with $u^N$.   

The order $N$ evaluation of the extended Yangian algebra $\mathcal{Y}_N$ is generated by the
coefficients appearing in the expansion of $T(u)$ \cite{AM07}. The expansion in $u$ of
$T(u+\delta)$ results in the equivalent algebra.

In general, the center of the $g\ell(n)$ type extended Yangian is generated
by the expansion in $u$ of the quantum determinant
\be \label{qdet}
 qdet T(u) = \sum_{\sigma} (-1)^{\sigma} T(u-n+1)_{1 \sigma(1)} ...
T(u)_{n\sigma(n)} = \ee $$
\sum_{\sigma} (-1)^{\sigma} T(u)_{\sigma(1) 1} ...
T(u-n+1)_{\sigma(n) n}.
$$
We have the simple co-multiplication property
\be \label{qdetT1T2}
 qdet T^1(u_1) T^2(u_2) = qdet T^1(u_1) \ qdet T^2(u_2).
\ee

We consider highest weight representations built on a singular vector $|0> $
obeying 
\be \label{hw<}
T_{ab}(u) |0> =0 , a < b, \ \ T_{aa}(u) |0> = \lambda_a(u) |0>.
\ee 
The related Verma module is obtained from $\mathcal{Y}_N$ as $\mathcal{Y}_N /
J_{\lambda}$, where $J_{\lambda}$ is generated by the expansion of $T_{ab}(u)  , a < b
$
and $T_{aa}(u) - \lambda_a(u) $.

The weight functions, polynomials in $u$ of order $N$, characterize the representation.
In this representation the quantum determinant appears as a number depending
on the weight functions. This number can be calculated as the eigenvalue of
the action on $|0>$. (\ref{qdet}) implies
\be \label{qdethw}
qdet T(u) |0> = \lambda_1(u) \lambda_2(u-1) ...\lambda_n(u-n+1) |0> .
\ee

In the case of first order evaluation
$ T(u) \to L(u) = I u + L$ we have $ \lambda_a(u) = u+ \lambda_a$, where
$\lambda_a$ are the $g\ell(n)$ representation weights and for the quantum
determinant we have 
\be \label{qdethwL}
 qdet L(u) |0> = (u+\lambda_1) (u-1+\lambda_2)  ... (u-n+1 + \lambda_n) |0>. \ee

In the order $N$ evaluation, where the 
 matrix of Yangian generators is composed as in (\ref{monodromy}),  
 the weight functions are obtained from the ones of the factors as
\be \label{lambdaau}
 \lambda_a(u) = \prod_{I=1}^N \lambda_a^I(u+\delta^I), \ \ \lambda^I_a(u) = u+
\lambda^I_a. 
\ee

\section{Jordan-Schwinger presentation of $g\ell(n)$}
\setcounter{equation}{0}

Using $n$ Heisenberg canonical pairs $x_a, \dd_a , a=1, ..., n$,  
\be \label{canon}
 [\dd_a, x_b] = \delta_{ab}, \ee
we build the Jordan-Schwinger (JS) presentation of the $g\ell(n)$ generators
\be \label{LJS+-}
 L^-_{ab} = -\dd_a x_b \ \ \ \ {\rm or } \ \ \ \ L^+_{ab} =  x_a \dd_b, \ee
both obeying the Lie algebra relations of $g\ell(n)$ and 
$ L^{\pm} (u) = I u + L^{\pm} $ obey the fundamental Yang-Baxter relation
(\ref{fundYB}).

In the following we shall prefer the version $-$ of (\ref{LJS+-})
and omit the superscript  in this case.

The center of the JS presentation of the $g\ell(n)$ algebra is generated by
\be \label{(xdd)}
(\underline x \underline \dd) = \sum_{a=1}^n x_a \dd_a. \ee
The center of $\mathcal{Y}_N$ (\ref{monodromy})  composed of JS factors
 is a subalgebra of the commutative algebra generated
by $ (\underline x^I \underline \dd^I), I=1, ..., N$.

The above statement about the center in the JS presentation is confirmed
by substituting (\ref{LJS+-}) to (\ref{qdet}). Indeed, we have
for $ a_1 < ... < a_m$
\be \label{da1x1}
 \sum_{\sigma} (-1)^{\sigma} \dd_{a_1} x_{\sigma(a_1)} ... \dd_{a_m}
x_{\sigma(a_m)} =  \dd_{a_1} x_{a_1} .
\ee
This implies that the result does not involve higher powers of $\dd_a$ or
$x_a$. Since the quantum determinant  commutes with $L_{ab}$ it has the form
$$ qdet L(u) = a(u) ( \underline \dd\underline x ) + b(u). $$

Evaluating the first form in (\ref{qdet}) with the substitution
$T(u) \to I u $ for all factors we obtain $a(u)$ and with the substitution
$T(u) \to I u  $ for all factors but one and 
$T(u) \to \dd_n x_n $ for one we obtain the result
\be \label{qdetLJS}
qdet L(u) = (u-n+1) ...(u-1) ( u- (\underline \dd \underline x ) ), 
\ee
and for the monodromy with the $L$ substituted in the JS form
$$ qdet T_N(u)  =  \prod_1^N (u+\delta^I-n+1) ...(u+ \delta^I-1) 
(u+\delta^I -n- ( \underline x^I  \underline \dd^I
).
$$

We consider representations on homogeneous functions of $x_a$ 
and fix the degree of homogeneity to $2\ell$. 
Considering the highest weight  condition in the form (\ref{hw<}) applied to
$L(u)$ in the JS presentation (\ref{LJS+-}) version $-$ 
we find the highest weight vector and the weight functions in explicit form 
\be \label{0>JS}
|0;2\ell> = x_n^{2\ell}, \lambda^{JS}_a(u) = u-1, a= 1, ..., n-1, 
\lambda^{JS}_n(u) = u-1-2\ell.
\ee 
The JS representation $L^{JS}_{ab} \cdot |0;2\ell> $ is spanned by
$$ \prod_{a=1}^{n-1} (x\p_a)^{m_a} x_n^{2\ell}, x\p_a = \frac{x_a}{x_n} $$
The rising operators produce powers of $x_a, a<n$ on expense of the degree
at $a=n$

In the set of homogeneous function of degree $2\ell$ we have the lowest
weight
vector $|\hat 0;2\ell> $ obeying
\be \label{lw}
 L_{ab} |\hat 0;2\ell> = 0, a>b, \ \  L_{aa} |\hat 0;2\ell> = \hat \lambda^{JS}_a |\hat 0;2\ell>  \ee
$$ (\hat \lambda^{JS}_a ) = (-2\ell-1, 0, ..., 0) $$
This lowest weight vector is represented by
$$ |\hat 0;2\ell> = x_1^{2\ell} = (x_1\p)^{2\ell} x_n^{2\ell} $$
The representation generated by the JS presentation $L_{ab}$ on $ |\hat 0;2\ell>$ is spanned by
$$ \prod_{a=2}^n (\hat x_a)^{ \hat m_a} =
(x_1\p)^{2\ell} x_n^{2\ell} \cdot (x_1\p)^{-\hat m_n} \prod_2^{n-1}
(\frac{x_a\p}{x_1\p})^{\hat m_a}. $$
The lowering operators produce powers of $x_a, a>1$ on expense of the degree
at $a=1$. 

In the case of non-negative integer $2\ell$ the rising operation in the
first case and the lowering  in the last case ends if the "source" at 
$a=n$ or $a=1$ is "exhausted".   
The highest and the lowest weight representations overlap in this case. 
We have the finite-dimensional representation with the restriction on the 
powers
$$\sum_1^{n-1} m_a \le 2\ell, \ \ \ \sum_2^n \hat m_a \le 2\ell. $$

 In view of the next section we introduce also the cyclic modification of
the above highest weight condition, $ k=0, 1, ..., n-1$ where $k=0$ refers to the
previous version (\ref{hw<}),
\be \label{hw<k}
T_{ab}(u) |0>^k = 0 {\rm \ \ for} \ \ \ 1-k \le a < b \le n-k \ee
$$ T_{aa}(u) |0>^k = \lambda^k_a(u) |0>^k, \ a=1, ..., n. $$
In the JS case, $T_{ab}(u) \to L_{ab}(u) = u \delta_{ab}-\dd_a x_b $,
the highest weight  vector is represented by $ |0;2\ell>^k = x_{n-k}^{2\ell} $
 The weight functions are $\lambda^{(k)}(u)_a = u+\lambda^{(k)}_a$ and the sequence of weights
is now
\be \label{lambdaJSk}
\lambda^{(k)}_1 = -1, ..., \lambda^{(k)}_{n-k-1} = -1, \lambda^{(k)}_{n-k} = -1-2\ell,
\lambda^{(k)}_{n-k+1} = -1 ...\lambda^{(k)}_n = -1.
\ee

The representation modules are represented by polynomial functions of the
ratios $ \varphi (x_a^{\prime (k)}),  x_a^{\prime (k)} = \frac{x_a}{x_{n-k}}
a\not = n-k $.  We shall use the restriction of the L operator to this
representation defined by
\be \label{Lu-u}
L(u) x_{n-k}^{2\ell} \varphi = x_{n-k}^{2\ell} L^{(k)} (u^-,u) \varphi,
\ \ \ u^- = u-2\ell .\ee
The additional argument involves the value of $( \underline x  \underline \dd )$
fixed to  $2\ell$. 

Below we shall have the constraints involving $x_a=0, a= n, n-1, ..., n-k+1$ 
On functions of the remaining coordinates the highest weight condition holds
in  the
standard form obeyed by $x_{n-k}^{2\ell} $.  
The lowest weight  conditions is  fulfilled by $x_1^{2\ell}$.

\section{Yang-Baxter relations}
\setcounter{equation}{0}

We consider the Yang-Baxter relation
\be \label{RLL}
R_{12}(u-v) L^{1 }_{a c}(u) L{^2 }_{c b}(v) =  L^{1 }_{a c}(v) L^{2 }_{c b}(u)
R_{12}(u-v).  
\ee
The form with explicit indices is presented to emphasize the difference to
(\ref{fundYB}). $L^I(u)$ is constructed from $x^I_a, \dd^I_a$ as above
(\ref{LJS+-}),
where we prefer the $L^-$ form. $R_{12}(w)$ is not a matrix and depends on
the canonical pairs. We consider (\ref{RLL}) as the defining relation for
$R_{12}$ and prove 
\begin{prop}
The integral over a closed contour in the complex plane
\be \label{R12w}
R_{12}(w) = \int \frac{dc}{c^{1-w}} e^{-c ( \underline x^1  \underline \dd^2 )
}, \ \ \ 
( \underline x^1  \underline \dd^2 ) = \sum_{a=1}^n x_a^1 \dd_a^2,  \ee
obeys (\ref{RLL}). 
\end{prop}

Proof:

We read (\ref{RLL}) as the defining relation of $R_{12}$ and solve it 
using the ansatz
$$ R_{12}(u) = \int \mathrm{d}c\, \phi(c)  
e^{-c (\underline{x}^1  \underline{\dd}^2)}. 
$$
We calculate the action of the shift operator on the product of $L$
matrices:
$$ e^{-c (\underline{x}^1  \underline{\dd}^2)}  \underline{\dd}^1 = 
(\underline{\dd}^1 + c \underline{\dd}^2 ) 
e^{-c (\underline{x}^1  \underline{\dd}^2)},
\qquad e^{-c (\underline{x}^1  \underline{\dd}^2)} \underline{x}^2 =
(\underline{x}^2 - c \underline{x}^1) \ e^{-c (\underline{x}^1  \underline{\dd}^2)},
$$
$$  e^{-c (\underline{x}^1  \underline{\dd}^2) } (L_1 (u) L_2
(v))_{ab}   e^{c (\underline{x}^1  \underline{\dd}^2)}= 
\big(L_{1 ac} (u) - c \dd^2_{ad} x^1_{dc}  \big)
\big(L_{2 cb} (v) + c \dd^2_{cd} x^1_{db} \big)
 =
$$ $$
\big( L^1 (v) L^2 (u) \big)_{ab} 
+ \big( L^2 (0) - L^1 (0)  \big)_{ab} \{ [u-v - c
(\underline{x}^1 \cdot \underline{\dd}^2) ] 
- c [ u-v +1 - c (\underline{x}^1  \underline{\dd}^2) ] \}. 
$$
The condition that the contribution of the second term 
vanishes implies a differential equation
on $\phi(c) $, because
$$ \int \mathrm{d}c\, \phi(c)  \cdot c\cdot (\underline{x}^1
\underline{\dd}^2)
e^{-c (\underline{x}^1  \underline{\dd}^2)} = 
- \int \mathrm{d}c\, \phi(c) \cdot  c \cdot \dd_c\,
e^{-c (\underline{x}^1  \underline{\dd}^2)} = 
\int \mathrm{d}c\, \dd_c  (c \phi(c) )  \, 
e^{-c (\underline{x}^1  \underline{\dd}^2)}. $$ 
In the last step we have assumed that the integration by parts is done
without boundary terms as it holds for closed contours.
Thus the condition on $\phi$ is
$$ \dd_c ( c \phi(c)) +  (v-u) \phi(c) = 0, $$
and it is solved by
$$ \phi(c) = \frac{1}{c^{1 + v-u} }. $$
The condition of vanishing of the third term can be written as
$$ 0 = \dd_c( c^2 \phi(c)) + c (v-u-1)\phi(c) = 
c [ \dd_c( c \phi(c)) + (v-u)\phi(c) ].
$$
We see that it does not imply a further condition on $\phi(c)$.
Thus we have proved that 
(\ref{R12w}) 
obeys the Yang-Baxter relation provided the simple rule of integration by
parts.
If the latter rule is different then the indicated procedure leads to the
appropriate modification.
\qed

\begin{prop}
The operator defined by integration over a closed contour
\be \label{R21w}
  R_{21} (w) = \int \frac{dc}{c^{1-w } } e^{-c
(\underline{x}^2\underline{\dd}^1)}
\ee
 with the argument substituted as $ \hat w= v-u +  (\underline{x}^1
\underline{\dd}^1) - (\underline{x}^2  \underline{\dd}^2) $ 
commutes with $L^1(u) L^2(v)$,
$$  R_{21}(\hat w)L^1(u) L^2(v) = L^1(u) L^2(v)  R_{21}(\hat w). $$
\end{prop}

Proof:
Similar to the previous proof we start with the calculation of the shift
operator action,
\be \label{ex2d1}
 e^{-c (\underline{x}^2\underline{\dd}^1)} L_{ac}^1(u) L_{cb}^2(v) =
(L_{ac}^1(u) +c \dd^1_ax_c^2) (L_{cb}^2(v) + c \dd^1_c x^2_b) e^{-c
(\underline{x}^2\underline{\dd}^1)}
= \ee 
$$
L_{ac}^1(u) L_{cb}^2(v) + c[ v-u+1 + (\underline{x}^1\underline{\dd}^1)
-(\underline{x}^2\underline{\dd}^2) + c\dd_c] \dd^1_a x^2_b e^{-c
(\underline{x}^2\underline{\dd}^1)}.
$$
We consider the ansatz
$$   R_{21} = \int dc \phi(c) e^{-c (\underline{x}^2\underline{\dd}^1)}
$$
and calculate 
$$   R_{21} L_{ac}^1(u) L_{cb}^2(v) = L_{ac}^1(u) L_{cb}^2(v)  \tilde R_{21}
+ \int dc c \phi(c) [ v-u+1 + (\underline{x}^1\underline{\dd}^1)
-(\underline{x}^2\underline{\dd}^2) + c\dd_c] \dd^1_a x^2_b e^{-c
(\underline{x}^2\underline{\dd}^1)}. 
$$ 
Integration by parts leads to the following condition of the vanishing of
the remainder 
$$ \dd_c (c \phi(c)) = [ v-u+1 + (\underline{x}^1\underline{\dd}^1)
-(\underline{x}^2\underline{\dd}^2) + c\dd_c] c\phi(c). $$ 
This fixes $\phi(c) $ in the form as asserted. 

\qed

\begin{prop}
The  $R$ operators (\ref{R12w}) and (\ref{R12w}) obey $RLL$ relations with the 
$L$ restricted to the
homogeneous functions  (\ref{Lu-u}).

\be \label{RLLu+u}
R_{12}(u-v) L^{1 (k) \ a}_c(u^-,u) L^{2 (l) \ c}_b(v^-, v) =  
L^{1 (k) \ a}_c (u^-, v) L^{2 (l) \ c}_b(v^-, u)
R_{12}(u-v),  \ee 
\be \label{RLLv+u}
R_{21}(u^--v^-) L^{1 (k) \ a}_c(u^-,u) L^{2 (l) \ c}_b(v^-, v) =  
L^{1 (k) \ a}_c(v^-, u) L^{2 (l) \ c}_b(u^-, v)
R_{21}(u^--v^-).  
\ee
\end{prop}

Proof:

We proof the second relation relying on Proposition 2. The first can be
proven in the analogous way.

We write the projected L operators as  
\be \label{LPi}
 L^{(k)} (u^-, u) = \Pi^{(k)} (2\ell) L(u) \Pi^{(k)} (2\ell), \ee
where the projection operator restricts the value of 
$(\underline x \underline \dd)$ to $2\ell$. This projection does not commute
with $R_{21}(w)$. Indeed, 
$$  R_{21}(w) (\underline x^1 \underline\dd^1) = \{ (\underline x^1 \underline \dd^1) -w \}  R_{21}(w), \ \ \  
  R_{21}(w) (\underline x^2 \underline \dd^2) = \{ (\underline x^2  \underline \dd^2) + w \}
R_{21}(w). $$
Therefore,
\be \label{R21Pi}
  R_{21}(w) \Pi^{(1, k)}(2\ell_1) \Pi^{(2, l)}(2\ell_2) =
\Pi^{(1, k)}(2\ell_1+w) \Pi^{(2, l)}(2\ell_2-w)  R_{21}(w). \ee

We obtain by (\ref{ex2d1}) and (\ref{R21Pi}) 
\be \label{R21L1L2}
  R_{21}(w)  L^{1(k)} (u^-, u) L^{2(l)} (v^-, v) =
 L^{1(k)} (u^-+w, u) L^{2(l)} (v^--w, v)  R_{21}(w) + \ee $$
\Pi^{1 (k)} (2\ell_1 +w)\Pi^{2(k)} (2\ell_2-w)
\int \frac{dc}{c^{-w}}  B_{ab}  e^{-(x^2 \dd^1)} \Pi^{1(k)} (2\ell_1)\Pi^{2(k)}
(2\ell_2),$$
$$  B_{ab} = 
[u-v- 1- (\underline x^1 \underline \dd^1)+
(\underline x^2 \underline \dd^2)+ c \dd_c]
\dd^1_a x^2_b.  
$$
By partial integration we obtain the analogous expression with $B_{ab}$
replaced by
$$\tilde B_{ab} = [u-v-  (\underline x^1 \underline \dd^1)+
(\underline x^2 \underline \dd^2)+ w]
\dd^1_a x^2_b.
$$
 The projectors allow to replace $\tilde B_{ab} $
by
$$ [u-v -2\ell_1 + 2\ell_2 +w ] \dd^1_a x^2_b. $$
Now we fix the value of $w$ to achieve the vanishing of the remainder,
$w= - u^-+v^-$. With this value of $w$ the first term on the r.h.s. of
(\ref{R21L1L2})
coincides with the r.h.s. of the Yang-Baxter relation (\ref{RLLv+u}). 

\qed


Consider the second order evaluation of the Yangian algebra, constructed by
the monodromy of two factors of L matrices (\ref{Lu-u})
$$ T^{(k,l)} (u; 2\ell_1,2\ell_2) = L^{1 (k)}(u-2\ell_1, u) L^{2 (l)}(u-2\ell_2, u)    
$$ 
with the weights 
$$ \lambda^1_a = -1, a\not = n-k, \lambda^1_{n-k} = -1-2\ell_1, $$ $$
\lambda_a^2 = -1, a \not = n-l, \lambda^1_{n-l} = -1-2\ell_2,
$$
and the one where the parameters are replaced by $2\ell_1\p ,2\ell_2\p   $.

\begin{prop}
A homomorphism between the algebras generated by $  T^{(k,l)} (u; 2\ell_1,2\ell_2)
$ and  $T^{(k,l) } (u; 2\ell_1\p,2\ell_2\p) $ exists, besides of the trivial
case $2\ell_i = 2 \ell_i\p + \delta $, only in the case of permutation of the
parameters, $2\ell_1\p= 2\ell_2, 2\ell_2\p = 2\ell_1$, where it is
represented by the operator (\ref{R21w}) $R_{21}(2\ell_1-2\ell_2) $.  

\end{prop}

{\it Proof:}

By the proposition 3 $R_{21}$ provides the wanted intertwining operator,
 $$ R_{21}(2\ell_1-2\ell_2) T^{(k,l)} (u; 2\ell_1,2\ell_2) = T^{(k,l)} (u;
2\ell_2,2\ell_1) R_{21}(2\ell_1-2\ell_2). $$

The ansatz of a product of $R_{12} (w\p)$ and $R_{21}(w)$ is sufficiently
general for solving
$$\mathcal {O} T^{(k,l)} (u; 2\ell_1,2\ell_2) = T^{(k,l)}(u;
2\ell_1\p,2\ell_2\p) \mathcal {O}. $$
$R_{12} $ acts trivially with $w\p =0$. In the previous proofs  we have seen how the
 R operator
arguments are fixed uniquely, in particular in the proof of proposition 3. 
\qed

\section{Biedenharn construction}
\setcounter{equation}{0}

 $L$ operators covering more general $g\ell(n)$ 
representations can be constructed from the monodromies (\ref{monodromy}) of $N
\le n$ factors imposing constraints. 
Let us describe in detail the first step of this construction.
\be \label{L1L2}
( L^1(u) L^2(u) )_{ab} = u^2\delta_{ab} - u ( \dd^1_a x^1_b + \dd^2_a
x^2_b ) + \dd^1_a ( \underline x^1  \underline \dd^2 ) x^2_b .
\ee
With the constraints 
$$ \mathcal{C}_{ab} = \dd^1_a ( \underline x^1  \underline \dd^2 ) x^2_b = 
( L^1(0) L^2(0) )_{ab} = 0
$$ 
we obtain the first order evaluation $L$ operator
\be \label{L1+2}
 L^1(u) L^2(u)|_{C_{ab}} = u L^{1+2}(u). \ee
The constraints are to be treated as  first order constraints in the sense of Dirac's
Hamiltonian formalism.
We solve the constraints by imposing 
\be \label{C2}
C^2 = ( \underline x^1  \underline \dd^2 ) = 0 \ee
The appropriate ordering rule is to be obeyed in the following, keeping
$\dd^1_a$ to the left of $( \underline x^1  \underline \dd^2 )  $ and 
$x^2_b$ to the right. 
We supplement (\ref{C2}) by the gauge condition
\be \label{gauge}
x^2_n = 0
\ee  
and express $\dd^2_n$ in terms of the other elements,
\be \label{d2n}
\dd^2_n= - (x^1_n)^{-1} \sum_{a=1}^{n-1} x^1_a \dd^2_a = - (x^1_n)^{-1}
( \underline x^1  \underline \dd^2 )^{(n)}. \ee
We use the notation
\be \label{()^k}
( \underline x^1  \underline \dd^2 )^{(k)} = \sum_{a=1}^{k-1} x_a^1 \dd_a^2
\ee
$L^{1+2}(u) $ has the first order evaluation form, 
$L^{1+2}(u) = uI + L^{1+2}$, and the explicit  matrix elements
are ($a,b, = 1, ..., n-1$),
\be \label{L1+20}
 L^{1+2}_{a,n} = -\dd^1_a x^1_n, \ \ 
  L^{1+2}_{n,b} = -\dd^1_n x^1_b + 
 (x^1_n)^{-1}( \underline x^1  \underline \dd^2 )^{(n)}x^2_b, 
\ee  $$ L^{1+2}_{n,n} = -\dd^1_n x^1_n, \ \ 
L^{1+2}_{a,b,} =  -\dd^1_a x^1_b - \dd^2_a x^2_b.
$$
$L^{1+2}_{ab} $ obey the $g\ell(n)$ Lie algebra relations.

The ordering means to keep $\dd^2_a$ to the left of $x^2_a$ 
before the constraint is substituted. In this way the initial constraints 
$(L^1(0) L^2(0))_{ab}= 0 $ holds also after the substitutions.  

In order to demonstrate the importance of the ordering we mention that the
ordering rule
changes, if we 
try to construct the first order evaluation starting from 

$$ (L^1(u-1) L^2(u))_{ab} = \delta_{ab} u (u-1) - u \dd^1_a x^1_b + (u+1)\dd^2_a x^2_b
 + \dd^1_a (\underline x^1 \underline \dd^2) x^2_b = $$ $$
(u-1) \left (\delta_{ab} u + L^1(0)_{ab} + L^2(0)_{ab}\right ) +
{C}_{r \ ab} $$

\be \label{Cabr}
C_{r \  ab} =\dd^1_a (\underline x^1 \underline \dd^2) x^2_b -
\dd^1_a x^1_b  = 
  \dd^1_a x^2_b (\underline x^1 \underline \dd^2) =
(L^1(0) L^2(1))_{ab} \ee
This way of construction works equally well with essential changes in the
details. 

Now we consider the restriction of the action on homogeneous functions.
The generators $L_{ab}^{1+2} $ obey the highest weight condition in the form
(\ref{hw<}) with
 $|0;2\ell_n,2\ell_{n-1}>^{1+2} $, the tensor product of
the h.w. vector $|0>^{1 (0)} $ obeying (\ref{hw<k}) with $k=0$ and 
$|0>^{2 (1)}$ with $k=1$,
\be \label{01+2} 
|0, 2\ell_n, 2\ell_{n-1} >^{1+2} = (x^1_n)^{2\ell_n} (x^2_{n-1})^{2\ell_{n-1}}, \ee
because
 the constraint $C^2$ (\ref{C2}), (\ref{gauge}) modifies the representation on functions of $x^I_a$ to
the ones  independent of $x^2_n$. 

 The sequence of weight functions   has the form $\lambda_a(u) = u+\lambda_a $, 
\be \label{weights1+2}
\lambda_1 = -2, ..., \lambda_{n-2} = -2, \lambda_{n-1} = -2-2\ell_{n-1},
\lambda_n = -1- 2\ell_n. \ee

We consider the L operator restricted to the action of homogeneous functions 
independent of $x^2_n$ 
of degree of homogeneity $2\ell_n$ and $2\ell_{n-1}$  and define in analogy
to (\ref{Lu-u})
the restricted L operator by
\be \label{L1+2lnln-1}
L^{1+2}(u) (x_n^1)^{2\ell_n} (x^2_{n-1})^{2\ell_{n-1}} \varphi =  (x^2_n)^{2\ell_n}
(x^2_{n-1})^{2\ell_{n-1}}
\mathcal{L}^{1+2} (u; 2\ell_{n}, 2\ell_{n-1} ) \varphi \ee
$ \varphi $ is a  function of the normal coordinates defined by the  ratios 
$ x_a^{\prime 1} = \frac{x^1_a}{x^1_n}, a=1, ..., n-1 $ and
$ x_a^{\prime 2 } = \frac{x^2_a}{x^2_{n-1}}, a=1, ..., n-2 $.
$\mathcal{L}^{1+2}$ can be expressed in terms of the restricted 
canonical pairs $x^{1 \prime}_1, \dd^{1 \prime}_a, a= 1, ..., n-1$,
and $x^{2 \prime}_a, \dd^{2 \prime}_a, a= 1, ..., n-2 $ of the normal
coordinates.

The construction can be continued stepwise to arrive at
$\mathcal{L} (u; 2\ell_n, ...,2\ell_1  )$,

\be \label{L1...Ln}
  \mathcal{L} (u; 2\ell_n, ...,2\ell_1  )
 = L^1(u+2\ell_n, u) L^2(u+2\ell_{n-1}, u)|_{C^2} ... L^n(u+2\ell_1,
u)|_{C^n}  
\ee
Here the action is restricted to   functions
$\varphi$ depending on $\half n (n-1)$ normal coordinates
\be \label{ratiok}
x_a^{1+k \prime} = \frac{x_a^{1+k}}{x_{n-i+1}^{1+k}}, 1\le a < n+1-k.\ee
 The sequence of weights is now complete, i.e. generic highest weight
representations can be described.
\be \label{weights1+...n}
\lambda_1 = -n-2\ell_1,  \lambda_{2} = -n+1-2\ell_2, ..., \lambda_{n-1} =
-2-2\ell_{n-1},\lambda_n = -1- 2\lambda_n. \ee
The shifts relating the weight $\lambda_a$ components to the homogeneity parameters
$2\ell_a$ are related to the components of the  half-sum of positive roots $\rho_a$
\be \label{rhoa}
\lambda_a + \rho_a + \half (n+1) + 2\ell_a = 0 
\ee

The condition $C^{1+k}, k=1, ..., n-1$ consists of the $k$ constraints
\be \label{Ck}
( \underline x^1  \underline \dd^{1+k} ) =0, ..., ( \underline x^{k}  \underline
\dd^{1+k} )= 0.
\ee
supplemented by the gauge conditions
$$ x^{1+k}_n = 0, ..., x^{1+k}_{n-k+1} = 0. $$
The conditions eliminate the canonical pairs $x^{1+k}_a, \dd_a^{1+k}, a= n-k+1, ...,
n$.
 
Further, the restriction to homogeneous functions fixes also 
$  ( x_{n-k}^{1+k}   \dd_{n-k}^{1+k} ) = 2\ell_{n-k}, k=0, ..., n-1 $. The remaining
degrees of freedom are represented by the ratio pairs
$ x_a^{1+ k \prime}$ (\ref{ratiok}), $\dd_a^{1+k \prime} = x_a^{1+k }
\dd^{1+k}_a, k=0, ..., n-1, a = 1, ..., n-k-1 $. 

We describe how the explicit form of $\mathcal{L} $ in terms of these
restricted canonical pairs can be obtained by solving the constraints and
restricting to functions of the normal coordinates. 

Consider first the case of the condition $C^2$ (\ref{C2}).
$ ( \underline x^1  \underline \dd^{2} ) =0 $ is solved by
$ \dd^2_n = - ( \underline x^{1 \prime}  \underline \dd^{2} )^{(n)} $.
The restriction to the functions, where the action restricts to $ x^2_{n-1} \dd^2_{n-1} =
2\ell_{n-1} $ results in
\be \label{d2n1}
 \dd^2_n = - \frac{1}{x^2_{n-1}}  \left ( 
( (\underline x^{1 \prime} - \underline x^{2 \prime})  \underline \dd^{2 \prime} )^{(n-1)} 
+ x^{1 \prime}_{n-1} 2\ell_{n-1} \right ). \ee

In general we solve (\ref{Ck}) for $\dd^{1+k}_n, ..., \dd^{1+k}_{n-k+1}$
taking into account also the conditions $C^{1+l}, l<k$. Thus we have
$$ ( \underline x^{1}  \underline \dd^{1+k} ) = 0,
( \underline x^{2}  \underline \dd^{1+k} )^{(n)} = 0, ..., 
( \underline x^{k}  \underline \dd^{1+k} )^{(n-k+2)} = 0.
$$
We separate the terms involving the unconstrained pairs only,
$$ \sum_{a=n-k+1}^n x_a^1\dd^{1+k}_a = - ( \underline x^{1}  \underline
\dd^{1+k})^{(n-k+1)},   $$
$$  \sum_{a=n-k+1}^{n-1} x_a^2\dd^{1+k}_a = - ( \underline x^{2}  \underline
\dd^{1+k})^{(n-k+1)},   $$
$$ .......$$
\be \label{dn-k+1} 
x^k_{n-k+1} \dd^{1+k}_{n-k+1} =   - ( \underline x^{k}  \underline
\dd^{1+k})^{(n-k+1)}.   \ee     
This system of $k$ linear equations is easily solved by inversion of the
involved triangular matrix. In the solution $x^1, ..., x^k$ can be replaced by 
the ratio coordinates (\ref{ratiok}).

The last factor in the monodromy (\ref{L1...Ln}) $L^n(u)$ supplies the last
parameter $\ell_1$, but the constraints $C^2, ..., C^n$ remove its canonical
pairs $x^n_a, \dd^n_a, a=2, ...n$ and the projection on the particular
homogeneous functions
 fixes even $x_1^n \dd_1^n$ to the value $2\ell_1$. The resulting matrix
  $L^n$ has non-vanishing elements in the first column only,
$$L^n_{11}(0) \to -1-2\ell_1, L^n_{1 a} (0)\to -\dd^{n *}_a x^n_1 , a= 2,
...n. $$ 
$\dd^{n *}_a, a=2, ..., n  $ are the solutions of (\ref{dn-k+1})  for the case
$k=n-1$, where the r.h.s of the $l$th line in (\ref{dn-k+1})  is now $-x_1^l
\dd^n_1$. 

We summarize the result as
\begin{prop}
 The first order evaluation $L$ matrix of the presentation of $g\ell(n)$
covering generic representations
can be constructed embedded in the algebra of $n^2$ Heisenberg canonical pairs  by
building the monodromy 
\be \label{Lgen}
\mathcal{L}(u; 2\underline \ell) =
L^1(u+2\ell_n, u) L^2(u+2\ell_{n-1}, u) ... L^n(u+2\ell_1, u)|_{\{C\} },  
\ee
$ 2\underline \ell = 2\ell_n, 2\ell_{n-1}, ..., 2\ell_1 $,
and imposing the conditions $\{ C\} = C^2, ..., C^n $ (\ref{Ck}).
The factors in the monodromy are the projections of JS type
$L$ matrices (\ref{LJS+-}), 
$$  L_{ab}^k(u-2\ell_{n-k}, u) = \Pi^{(k)}(2\ell_{n-k+2}) 
\left ( u \delta_{ab} - \dd^k_a x_b^k \right )  \Pi^{(k)}(2\ell_{n-k+2}), $$
where the projection is described by (\ref{Lu-u}).
\end{prop}

An analogous construction  of a generic $g\ell(n)$ representation can be 
formulated working with the constraints in the form (\ref{Cabr}). In this case
the factors in the  modification of (\ref{Lgen}) appear with shifts in $u$.

\section{Evaluation of order $N$}
\setcounter{equation}{0}

\subsection{The parameters of the presentation}

We consider the order $N$ evaluation of the $g\ell(n)$ Yangian algebra
$\mathcal{Y}_N$ by
the presentation in terms of the monodromy built as the $N$-fold product of copies
of $\mathcal{L}(u; 2\underline {\ell} ) (\ref{Lgen}), 
2 \underline \ell=  (2\ell_n, ..., 2\ell_1 ) $,
\be \label{monogen}
T_N(u, 2\frak l) = \mathcal{L}^{1 }(u; 2\underline {\ell}^1 ) ... \mathcal{L}^{N }(u; 2\underline
{\ell}^N ).
\ee

It is embedded in the Heisenberg algebra of the $n^2 N$ canonical pairs 
$$ x_a^{(I,i)}, \dd_a^{(I,i)}, \ \ I=1, .., N; i=1, .., n; a=1, .., n. 
$$
The constraints $ C^{(I,i)} $ eliminate the pairs with $ a>n-i+1$. 
The factors $\mathcal{L}^{I }$ of the monodromy are the L matrices of generic 
$g\ell(n) $ representations constructed above. 
In particular the representations have been restricted to
homogeneous functions of the form 
$$ \prod_{i=1}^n (x_{n-i+1}^{(I,i)})^{2\ell^I_{n-i+1}}\varphi, $$
where $\varphi$ depends on the normal coordinate ratios 
\be \label{ratio}
x_a^{(I,i) \prime} = \frac{x_a^{(I,i)}}{x_{n-i+1}^{(I,i)} },  a= 1, ...,
n-i.
\ee
This reduction allows the expression in terms of the corresponding reduced
canonical pairs $ x_a^{(I,i) \prime }, \dd_a^{(I,i) \prime }, 
a= 1, .., n-i $. 
The constant function $1$ obeys the highest weight conditions in the form
(\ref{hw<}) and 
the weight functions  are
\be \label{weightf}
 \lambda_a(u, 2\frak l) = \prod_{I=1}^N \lambda_a^I(u), \ \ \lambda_a^I(u)
= u+n-a+1 +2\ell_a^I. \ee
 
In this way the  ordered $nN$ parameter array
$$ 2 \frak l = 2\underline {\ell}^1; ...;  2\underline {\ell}^N = 2\ell^1_1 , ...,
2\ell^1_n; ...; 2\ell^N_1, ..., 2\ell^N_n $$
characterizes the  presentation of the algebra $\mathcal{Y}_N(g\ell(n),
2\frak l)$. 

We calculate the quantum determinant (\ref{qdet}). 
$$ qdet T_N(u, 2 \frak l) = \prod_I qdet \mathcal{L}^{I }(u; 2\underline
{\ell}^I )
$$
By (\ref{qdethwL}) and (\ref{weightf})
\be \label{qdetT2l}
 qdet \mathcal{L}^{I }(u; 2\underline {\ell}^I ) \cdot 1
= \prod_a (u+n+2\ell_a^I) .
\ee

The result is symmetric  in the representation
parameters $2\ell_i^I $, the values of $(\underline{x}^(I,i) \underline
\dd^{I,i})$ in action to the homogeneous functions.     

\subsection{Permutations by $R$ operators}

We discuss, how the R operators constructed above allow to obtain operators 
acting on the monodromy matrices  as
$$ \mathcal{S}  T_N(u, 2 \frak l) =  T_N(u, \sigma 2 \frak l)  \mathcal{S}
$$
where $\sigma 2 \frak l$ denotes a permutation of the parameter array 
$2 \frak l$. 

\begin{prop}
The Yang-Baxter operator  (\ref{R21w})
provides by appropriate substitutions 
operators of adjacent parameter permutations, intertwining the corresponding
algebra presentations. 
The operator $S^I_i$
permuting adjacent parameters $2\ell^I_{n-i+1}, 2\ell^I_{n-i}$
in $\mathcal{L}^{I }$
is obtained from
\be \label{RIi+1}
 R_{(I,i+1),(I,i)}(2\ell^I_{n-i+1}-2\ell^I_{n-i}) = \int 
\frac{dc}{c^{1-2\ell^I_{n-i+1}+ 2\ell^I_{n-i}}} 
e^{-c(\underline{x}^{I,i+1} \underline{\dd}^{I,i})}
\ee
by imposing the constraints $C^{I,i+1}, C^{I,i}$ on the on the involved
canonical pairs.
The operator $S^I$ permuting adjacent parameters 
$2\ell^I_{1}, 2\ell^{I+1}_{n}$
of adjacent $\mathcal{L}^{I }$ and $\mathcal{L}^{I+1}$
is obtained from
\be \label{RI1I+1}
R_{(I+1,1,(I,n)}(2\ell^I_{1} -2\ell^{I+1}_{n}) = 
\int \frac{dc}{c^{1-2\ell^I_{1}+2\ell^{I+1}_{n}}} 
e^{-c(\underline{x}^{I,i+1} \underline{\dd}^{I,i})}
\ee
by imposing the constraints $C^{I,2}, ..., C^{I,n}$  on the involved
canonical pairs.
\end{prop}

Proof:

The considered monodromy has been constructed by ordered matrix products of 
$n N$ factors of the JS type L matrices $L(u^-, u)$,
$$ T^0_N (u) =  T^1_n (u) ... T^N_n(u), \ \ \
T^I_n(u) = L^1(u-2\ell^I_n, u) ... L^n (u-2\ell^I_1, u) $$
by imposing the constraints (\ref{Ck}),
$$ T_N(u) = T^0_N(u)|_{\{C^{(I,i)}\} }. $$
Before the constraints are imposed the above 
$ R_{(I,i+1,(I,i)}( 2\ell^I_{n-i}-2\ell^I_{n-i+1})$
 interchanges $2\ell^I_{n-i+1}, 2\ell^I_{n-i}$ and 
$ R_{(I+1,1,(I,n)}(2\ell^{I+1}_{n}-2\ell^I_{1})$ interchanges 
$2\ell^I_{1}, 2\ell^{I+1}_{n}$ as established in the propositions 2 and 3.

The constraints affect the argument of the involved exponential.
For consistency it is important that the constraints commute with the
corresponding $R$ operator. To see this, it is sufficient to consider a particular case.
Similar as in the proof of (\ref{RLLv+u}) we have
$$ R_{21}(w) \Pi^1(2\ell_1 ) \Pi^2(2\ell_2) (\underline x^1 \underline \dd^2) = 
\Pi^1(2\ell_1+w) \Pi^2(2\ell_2)-w) (\underline x^1 \underline \dd^2)  R_{21}(w) + $$ $$
\Pi^1(2\ell_1+w) \Pi^2(2\ell_2 -w)
\int \frac{dc}{c^{-w}} \{ (\underline x^1\underline \dd^1) - 
(\underline x^2\underline \dd^2)  - w \}
e^{-(\underline x^2 \underline \dd^1)} $$
The remainder vanishes indeed at $w=2\ell_1-2\ell_2$.
This confirms that the constraints do not change.  Thus the  parameter permutations
in the monodromy $T_N(u, 2\frak l) $
are caused by
the action of the corresponding $R$. 
 \qed

We comment on the parameter permutations in the case of first order
evaluation, $N=1$.
The weights  of the  JS representation considered in sect. 2 have the
particular form (\ref{0>JS})
$$ ( \lambda_a )^{JS} = ( -1,-1, ...,-1,-1-2\ell) $$
Compared to the generic form constructed in sect. 5,
$$ ( \lambda_a )^{B} = ( -n-2\ell_1,-n+1-2\ell_2, ..., -2-2\ell_{n-1}, -1-2\ell_n ) $$
this corresponds to the particular sequence of representation parameters 
$$(2\ell_a)^{JS} = (-n+1,-n+2, ...,-1, 2\ell). $$
Other representations with one generic parameter $2\ell$
obtained by the parameter permutation operations cannot in general presented
in the same JS form. The particular one with the parameter  $2\ell$ shifted cyclically
to the first position
$$(2\ell_a)^{cycl} = (2\ell, -n+1, ...,-1 ) $$
has the weight sequence
$$ ( \lambda_a )^{cycl} = ( -n-2\ell, 0, ..., 0, 0) $$. It 
can be presented in the JS $L^+$ version of (\ref{LJS+-}) with the same form of the h.w.
condition (\ref{hw<}). In this case the highest weight vector is represented by $x_1^{2\tilde
\ell}$ and
the weights are
$$ (\lambda_a)^{JS+} = (2\tilde \ell, 0, ..., 0). $$
Thus the representation obtained by the cyclic permutation $R$ operations from the $L^-$
JS representation is equivalent to the $L^+$ JS representation if
the parameters are related as $2\tilde \ell = -2\ell-n$.


\subsection{Operators of parameter permutations}

Let us consider  simple cases of parameter permutations 
in order to provide explicit expressions. 

In the case of the first step of the above construction (\ref{L1+2lnln-1}),
$$ S_{12} \mathcal{L}^{1+2}(u;2\ell_n, 2\ell_{n-1}) =  \mathcal{L}^{1+2}(u;2\ell_{n-1}, 2\ell_{n})
S_{12},$$ 
the intertwiner with the representation with the permutation of $2\ell_{n},
2\ell_{n-1}$ is obtained as
 $ R_{21}( 2\ell_{n}-2\ell_{n-1} )|_{C^2} $. 
The gauge condition in the constraint changes the argument in the
exponential in $R_{21}$ as
\be \label{S12}
 (\underline{x}^2 \underline{\dd}^1)|_{C^2} =
 (\underline{x}^2 \underline{\dd}^1)^{(n)} = 
 \frac{x^2_{n-1}}{x^1_n} (\underline{x}^{2 \prime} \underline{\dd}^{1
\prime})^{(n)} = 
 \frac{x^2_{n-1}}{x^1_n} ((\underline{x}^{2 \prime} \underline{\dd}^{1
\prime})^{(n-1)} + \dd^{1 \prime}_{n-1} ) \ee 

Consider next a monodromy of such a Biedenharn first step factor and a unconstrained JS
factor (\ref{Lu-u}),
$$ S_{23} \mathcal{L}^{1+2}(u;2\ell_{n}, 2\ell_{n-1}) L^3(u-2\ell_3, u) =
 \mathcal{L}^{1+2}(u;2\ell_n, 2\ell_3) L^3(u-2\ell_{n-1}, u) S_{23}. $$
This intertwiner $S_{23}$  permutating  $2\ell_{n-1},
2\ell_3 $ is obtained as
$ R_{32}(  2\ell_{n-1}- 2\ell_{3})|_{C^2} $.

The gauge condition changes the argument in the exponential as
$$ (\underline{x}^3 \underline{\dd}^2)|_{C^2} = 
x^3_n \dd^2_n|_{C^2} x^3_{n-1} \dd^2_{n-1} + 
(\underline{x}^3 \underline{\dd}^2)^{(n-1)} = $$ $$
- \frac{x^3_n}{x^2_{n-1}} \left ( 2\ell_{n-1} x^{1 \prime}_{n-1} +
((\underline{x}^{1 \prime}-\underline{x}^{2 \prime}) \underline{\dd}^{2 \prime})^{(n-1)} 
\right ) + $$ $$ \frac{x^3_{n-1}}{x^2_{n-1}} (2\ell_{n-1}-(\underline{x}^{2 \prime}
\underline{\dd}^{2 \prime})^{(n-1)}) + 
 \frac{x^3_{n}}{x^2_{n-1}} 
(\underline{x}^{3 \prime} \underline{\dd}^{2 \prime})^{(n-1)} = $$ 
\be \label{S23}
\frac{x^3_{n}}{x^2_{n-1}} \left ( 2\ell_{n-1} ( x^{3 \prime}_{n-1} - x^{1
\prime}_{n-1} ) - (\underline{x}^{1 \prime} \underline{\dd}^{2 \prime})^{(n-1)}
+ (\underline{x}^{3 \prime} \underline{\dd}^{2 \prime})^{(n-1)} 
\right ).  
\ee

Consider now the general case of a monodromy (\ref{monogen}).
  
Let us find  the explicit expressions of $S^I_i$ for
intertwining the  presentations of the algebras with the permutation of
$2\ell^I_{n-i+1}$ and
$2\ell^I_{n-i}$, 
$  R_{(I, i+1), (I,i)}|_{C^{(I,i)}, C^{(I,i+1)}} $ (\ref{RIi+1}), 
by calculating the exponent in the involved shift operator 
under the constraints. Only the restrictions on the $x$ components 
$ x_n^{I,i+1} = 0,x_{n-1}^{I,i+1}= 0 , ..., x_{n-i+1}^{I,i+1}= 0  $
matter now, 
$$ (\underline{x}^{I,i+1} \underline{\dd}^{I,i})|_{C^{(I,i)}, C^{(I,i+1)}} =
 (\underline{x}^{I,i+1} \underline{\dd}^{I,i})^{(n-i+1)} = $$ 
\be \label{Si,i+1}
\frac{x^{I,i+1}_{n-i}}{x^{I,i}_{n-i+1} } \left 
( (\underline x^{I, i+1 \prime} \underline \dd^{I,i} )^{(n-i)} + \dd^{I,i
\prime}_{n-i} \right ), \ee
We have introduced the ratio canonical pairs (\ref{ratio}). 
This generalizes (\ref{S12}).
By a substitution of the
integration variable $c$ the factor $\frac{x^{I,i+1}_{n-i}}{x^{I,i}_{n-i+1} }$
can be removed from the argument of the shift operator exponential and
appears as the appropriate factor for changing the highest weight  vector
$$ \left (\frac{x^{I,i+1}_{n-i}}{x^{I,i}_{n-i+1} } \right )^{2\ell^I_{n-i+1} -
2\ell^I_{n-i} } \ |0; ...,2\ell^I_{n-i+1}, 2\ell^I_{n-i}, ...> =
|0; ..., 2\ell^I_{n-i}, 2\ell^I_{n-i+1}, ... >.
$$

Let us find also the explicit expressions of  $S^I$
for
intertwining the presentations of the algebras with the permutation of $2\ell^I_1$
and $2\ell^{I+1}_n$, (\ref{RI1I+1}).
No constraints are imposed on the pairs $(x_a^{I+1,1}, \dd_a^{I+1,1})$ 
but the pairs $(x_a^{I,n}, \dd_a^{I,n} $ are affected by all 
 $C^{I,2}, ..., C^{I,n}$.
 
$$  (\underline{x}^{I+1,1} \underline{\dd}^{I,n})|_{C^{I,2}, ..., C^{I,n}} =
\sum_{a_2}^{n} x_a^{I+1,1} [\dd^{I,n}_a ]_{C^{I,2}, ..., C^{I,n}}
+  x_1^{I+1,1} \dd_1^{I,n}  $$
$[\dd^{I,n}_a ]_{C^{I,2}, ..., C^{I,n}}$ are obtained by 
 (\ref{dn-k+1}) in the case $k=n-1 $. We have to solve
$$  (\underline{x}^{I,1} \underline{\dd}^{I,n}) = 0, \ 
 (\underline{x}^{I,2} \underline{\dd}^{I,n})^{(n)} = 0,\  
...,
 (\underline{x}^{I,n-1} \underline{\dd}^{I,n})^{(3)} = 0.
$$


This can be written in matrix form as
$$ 
\begin{pmatrix} 
x_{n}^{I,1} &  x_{n-1}^{I,1 } &...& x_2^{I,1 }  \\
0           &   x_{n-1}^{I,2 }&...& x_2^{I,2 } \\
0           & 0            &...& ...                \\
...         &    ...                 &...& ...                \\
0           & 0                      &...&  x_2^{I,n-2 } \\   
0           & 0                      &...& x_2^{I,n-1 }          \\
 \end{pmatrix}
\begin{pmatrix}
 \dd^{I, n}_n \\
\dd^{I, n}_{n-1}  \\
... \\
... \\
\dd_{3}^{I,n} \\
\dd_{2}^{I,n } \\ 
\end{pmatrix}
= 
- 
\begin{pmatrix}
 x^{I, 1}_1 \\
x^{I, 2}_{1}  \\
... \\
... \\
x_{1}^{I,n-2} \\
x_{1}^{I,n-1 } \\ 
\end{pmatrix}
\dd^{I n}_1 .
$$

The components $a=2, ..., n$ of $\underline{\dd}^{I,i}$ result in
expressions in
the ratio coordinates $x^{I,i \prime}$
$ i= 1, ..., n-1, a= 1, ..., n-i $, and $x_1^{I,n}, \dd_1^{I,n}$

We introduce the representation condition
$$ x_1^{I,n} \dd_1^{I,n} \to 2\ell_1^I $$
and obtain the result in coordinate components only,

\be \label{SI,I+1}  
(\underline{x}^{I+1,1} \underline{\dd}^{I,n})|_{C^{I,2}, ..., C^{I,n}} =
2\ell_1 \frac{x^{I+1,1}_n}{x_1^{I,n} } (x_1^{I+1,1 \prime}- X ),
\ee $$
X=
\begin{pmatrix} 
1 & x_{n-1}^{I+1,1 \prime}& ... & x_3^{I+1,1\prime} & x_2^{I+1,1\prime} 
\end{pmatrix}
{\begin{pmatrix} 
1           &  x_{n-1}^{I,1 \prime} &...& x_2^{I,1 \prime}  \\
0           &    1                   &...&  x_2^{I,2 \prime} \\
0           & 0                      &...& ...                \\
...         &    ...                 &...& ...                \\
0           & 0                      &...&  x_2^{I,n-2 \prime} \\   
0           & 0                      &...& 1                  \\
 \end{pmatrix} }^{-1}
\begin{pmatrix}
 x_{1}^{I,1 \prime} \\
x_{1}^{I,2 \prime}  \\
... \\
... \\
x_{1}^{I,n-2 \prime} \\
x_{1}^{I,n-1 \prime} \\ 
\end{pmatrix}.
$$
Again, the factor $2\ell_1 \frac{x^{I+1,1}_n}{x_1^{I,n} }$ can be removed
from the argument of the shift exponential by substitution in the integral,
resulting in the appropriate factor in front to change the 
highest weight vector. 

Removing this factor the remaining expression in the normal coordinates is 

\be \label{SIdet}
  (\underline{x}^{I+1,1} \underline{\dd}^{I,n})|_{C^{I,2}, ..., C^{I,n}} =
const \ \ det (\underline{x}^{I,1}
\underline{x}^{I,2}...\underline{x}^{I,n-1} \underline{x}^{I+1,1} )
\ee

This can be obtained by writing the matrix inverse in terms of minors.

In the following we restrict to the normal coordinates and the corresponding
canonical pairs.  In the following, simplifying the notation we omit the primes, which 
distinguished so far the normal from the homogeneous coordinates.

We introduce notations for the considered permutations of adjacent items in
the list of permutation parameters $2\frak l$ leaving all others unchanged:
\be \label{sigmaI}
 \sigma^I_i : 2\ell^I_{i+1} \leftrightarrow 2\ell^I_{i}, \ \ 
\sigma^I :  2\ell^I_1 \leftrightarrow 2\ell^{I+1}_{n}
 \ee

The obtained explicit forms (\ref{Si,i+1}), (\ref{SI,I+1}), (\ref{SIdet})  are summarized as

\begin{prop}
In action on the monodromy $T(u,2\mathfrak l)$
$$ T(u,2\sigma \mathfrak l) = S T(u,2\mathfrak l)  S^{-1} $$
the particular adjacent permutation $\sigma^I_i$ in $2\mathfrak l$ permuting 
$ 2\ell^I_{n-i+1}, 2\ell_{n-i}i ,  i=1, ...,n-1 $ is represented by
\be \label{SIi}
 S^I_{i,i+1}(2\ell^I_{n-i+1}-2\ell^I_{n-i}), 
 S^I_{i,i+1}(w)  = \int dc c^{w-1} \exp (-c(\underline x^{I,i+1} \underline
\dd^{I,i})) \ee
where the $n$ component vectors are 
$$ \underline{ x}^{I,i} = (x^I_1, ...,x^I_{n-i}, 1, 0, ...,0), \ \ \ 
\underline \dd^{I,i} = (\dd^I_1, ...,\dd^I_{n-i},0, 0, ...,0) $$
and the adjacent permutation $\sigma^I$ permuting $2\ell^I_1, 2\ell^{I+1}_n $
\be \label{SI}
 S^I(2\ell^I_1-2\ell^{I+1}_n), 
S^I(w) = \left (det (\underline x^{I,1}, \underline x^{I,2}, ..., \underline
x^{I,n-1}, \underline x^{I+1,1} ) \right )^w
\ee 
\end{prop}

The Yangian algebra relation known in general setting is now appearing 
in the particular presentation of the explicit Biedenharn construction
based on Heisenberg algebras.

\begin{prop}
The order $N$ evaluations of the Yangian algebra generated by 
$ T(u, 2\frak l ) $ as defined in (\ref{monogen}) characterized by the parameter set
$$
2\frak l = (\underline 2\ell^1, \underline 2\ell^2. ..., \underline
2\ell^N ), \ \ 2 \underline \ell^I = 2\ell^I_n, 2 \ell^I_{n-1}, ...,2\ell^I_1 $$
and  the representation characterized by the parameter set 
$ 2\frak l\p $ are related  by a coinciding center and by the existence of
morphisms 
$$ T(u, 2\frak l\p ) = \frak S (2\frak l, 2\frak l\p) T(u, 2\frak l )
\frak S (2\frak l\p , 2\frak l) $$
iff the parameter sets are related be permutation of the weight parameter
components up to an overall shift of the parameters 
$$ 2\frak l\p =  \sigma \ 2\frak l \to 2\frak l +
v . $$ 
\end{prop}

Proof:

The center is generated by the expansion in $u$ of the quantum determinant.
Its dependence on the parameters has been obtained in (\ref{qdetT2l}) and shown to be
symmetric with respect to parameter permutations. 

If the parameter sets are related by permutation we choose a decomposition
into elementary permutations of adjacent parameters
$2\ell^I_{k+1}, 2\ell^I_k, k= 1, ..., n-1$ or $2\ell^I_n, 2\ell^{I+1}_1 $
and construct $\frak S (2\frak l, 2\frak l\p)$ correspondingly by the 
permutation operators
$ {S}^I_k (2\ell^I_{k+1}- 2\ell^I_k), 
 {S}^I (2\ell^I_{1}- 2\ell^{I+1}_n) $
given explicitly above, proposition 6 and  (\ref{Si,i+1}), (\ref{SI,I+1}). 

These operators of elementary permutations have been obtained from the 
$R$ operators obeying the Yang-Baxter relations.   
If the parameter sets coincide besides of two adjacent items $2\ell,
2\ell\p$ and $2\tilde \ell,
2\tilde \ell\p$ then 
 the Yang-Baxter relation (\ref{RLLv+u}) holds for of the corresponding factors
$L^{JS} (u-2\ell,u) L^{\prime JS} (u-2\ell\p,u) $
and $ L^{JS} (u-2\tilde \ell,u) L^{\prime JS} (u-2 \tilde \ell\p,u) $
only if the parameters in the second  are permutations of the ones in the
first, according to proposition 4.
\qed

\section{Morphisms of representations } 
\setcounter{equation}{0}

The generators $T(u, 2\frak l)$ are constructed in the frame of the
Heisenberg algebra generated by $x^{Ii}_a, \dd^{Ii}_a, I=1, ..., N,
i=1, ..., n-1, a=1,..., n-i $. 
We consider the representations on functions of the normal coordinates
$ x^{Ii}_a $, where $x^{Ii}_a, \dd^{Ii}_a $ act as multiplication and
derivative respectively. The function $1$ constant in all arguments obeys the 
highest weight  condition (\ref{hw<})
\be \label{hw1}
 T_{ab}(u, 2\frak l) \cdot 1 = 0, a<b,  \ \ \ 
T_{aa}(u, 2\frak l) \cdot 1 = \lambda_a(u, 2\frak l), \ee
where the weight functions are given in (\ref{weightf}).

The result of the action  $T_{ab}(u, 2\frak l) $ on $1$ lies in the subspace
of functions spanned by the monomials in $ x^{I,i}_a $.

There are more functions of the normal coordinates obeying the highest
weight relation (\ref{hw<}). 
Consider along with the parameter sequence $2\frak l$ a further
one $\frak l\p = \sigma^I \frak l$, where the adjacent $2\ell^{I}_1, 2\ell^{I+1}_n $
are permuted. Then the function
\be \label{phiI}
 \phi^I =  S^{I}(2\ell^{I+1}_n- 2\ell^I_1 ) \cdot 1 = 
= \left (det (\underline x^{I,1}, \underline x^{I,2} ...\underline x^{I,n-1
, \underline x^{I+1,1}}) \right )^{ 2\ell^{I+1}_n - 2\ell^I_1} \ee
   obeys 
\be \label{hwphi}
 T_{ab}(u, 2\frak l) \cdot \phi^I = 0, a<b, \ \ \  
 T_{aa}(u, 2\frak l) \cdot \phi^I = \lambda_a(u, 2\frak l\p ) \phi^I \ee
This is obtained by the action of the permutation operator 
$ S^{I}(2\ell^{I+1}_n- 2\ell^I_1 ) $ on (\ref{hw1}) with
$ 2\frak l$ replaced by $2\frak l\p $ and using Proposition 7.

A further example is obtained
 by the action of the permutation operator 
$ S^I_i$. We start with (\ref{hw1}), where $2\frak l$ is replaced by
$\sigma^{-1} 2 \frak l $ with $\sigma^{-1} $ chosen such that
 the part of the sequence $2\ell^I_2,2\ell^I_1; 2\ell^{I+1}_n $
is
permuted to $2\ell^1_1,2\ell^{I+1}_n; 2\ell^I_2 $. We act  with the permutation
operator 
$$ \frak S (\sigma^{-1} 2 \frak l) = S^I_1(2\ell^1_2-2\ell^1_1)
S^I(2\ell^{I+1}_n- 2\ell^1_2) $$
and obtain the highest weight relation of the form (\ref{hwphi})
with $\phi^I$ replaced by  the function
\be \label{phiI1}
 \phi^I_1 =
const \ \ 
\left (det (\underline x^{I,1}, \underline x^{I,2} ..., \underline x^{I,n-2}, \underline
x^{I,n}, \underline x^{I+1,1}) \right)^{2\ell^I_1- 2\ell^I_2} \ee 
$$
 (det (\underline x^{I,1}, \underline x^{I,2} ...\underline x^{I,n-2},\underline x^{I,n-1} , 
\underline x^{I+1,1} ) )^{2\ell^{I+1}_n - 2\ell^I_1} $$

$\underline x^{I,n} $ is the unit vector with all components zero besides the
first, actually not depending on $x$. 

The above examples are particular cases of the following scheme.
The presentation of the order $N$ evaluation of the Yangian
constructed depends of the parameter array 
$2\frak l = (2 \underline \ell^1, ..., 2 \underline \ell^N) $.
Let $\sigma $ be a permutation on this set.
Let $\frak S$ be a sequence of the elementary permutation operators
corresponding to this parameter permutation. 
Consider the highest weight  condition and the representation generated by
$T(u;2\sigma^{-1} \frak l) $ on $1$.
$$ T_{ab} (u;2\sigma^{-1} \frak l) \cdot 1 = 0, a<b,$$
$$ T_{aa} (u;2\sigma^{-1} \frak l) \cdot 1 = \lambda_a(u;2\sigma^{-1} \frak l) \cdot
1,  $$
The representation is generated by the action of 
$ T_{ab} (u;2\sigma^{-1} \frak l) , a>b $ on  $1$. 
Let us  act with $\frak S(\sigma^{-1} \frak l) $  
to obtain
$$ T_{ab} (u;2 \frak l)  \frak S  (\sigma^{-1} \frak l) \cdot 1 = 0, a<b, \ \ \ 
 T_{aa} (u; \frak l) \frak S(\sigma^{-1} \frak l) \cdot 1 = \lambda_a(u;2 \frak l)  \frak S(\sigma^{-1} \frak l) \cdot
1  $$
and obtain the representation generated by 
$ T_{ab} (u;2 \frak l) , a>b$ acting on $ \phi=\frak S(\sigma^{-1} \frak l)
\cdot 1$.  

Thus $\frak S $ maps the representation generated by $T(u;2\sigma^{-1} \frak l)$
on $1$ to the representation generated by $T(u;2 \frak l)$ on $\Phi = \frak S \cdot
1$. 

Obviously by this scheme  further examples of functions obeying the highest
weight condition can be
obtained. The number of such highest weight
functions depends on the rank, e.g. at $n=2$ the first determinant factor in
(\ref{phiI1}) is constant and thus $\phi^I$ coincides with $\phi^I_1$ in
this case. 

We see how 
the permutation operators can result in morphisms of the representations of
the algebra generated by $T(u, 2\frak l)$ spanned by the monomials
and characterized by the weight functions $\lambda(u, 2\frak l)$ and the
representation spanned by monomials multiplying $\phi(x, 2\frak l\p) $
characterized by the weight functions $\lambda(u, 2\frak l\p)$. 
For generic parameters we have highest weight representations, 
in $T_{ab}(u, 2\frak l) \cdot \phi$ there is only one highest weight vector
represented by $\phi$. The representation space includes all the monomials
multiplied by $\phi$. 
We notice the
interesting case, where the differences of the representation parameters
appearing as exponents in the expressions of $\phi$ are non-negative integers.  
Then $\phi$ lies in the subspace spanned by monomials and the
representation spaces overlap.  In this case we have both
 $1$ and $\phi$ as highest weight vectors. 

The situations where not all 
monomials appear in $T_{ab}(u, 2\frak l) \cdot \phi$ are related with 
non-negative integer parameter differences too. 

The constant appearing in the second permutation step (\ref{phiI1})
depends of the representation parameters. 
Zeros and poles appear just for integer differences.
Thus the analytic structure of this permutation coefficient
indicates the degeneracy in the representation. We shall focus on this
aspect in the following.


The specification of the action of $S^{I}_i $ (\ref{SIi}) on $1$ requires  careful 
limits. Avoiding  complications we restrict to sequences of 
permutations, called beta sequences, where the involved $ S^I_{i} $ with their shift
operators  
act non-trivially leading to well defined  Beta integrals.
 
If in particular 
$\frak S \cdot 1 = \pi(2\mathfrak l) 1  $ we call it a 1 to 1 beta sequence.
 In this case the representations generated by the action of $ T(u;2\sigma^{-1} \frak l)
$ on
$1$ is mapped to the representation generated by $ T(u; \frak l) $ on $1$.
 The permutation coefficient $\pi(2 \mathfrak l)$ is given by 
a product of Beta functions.

In the example of the permutation 
of $2\ell^1_n $ and $2\ell^{2}_n $
we have the sequence (\ref{sigmaI}) 
$$ \sigma^{12}_n = = \sigma^1 \sigma^1_{n-1} \sigma^1_{n-2} ... \sigma^1_{2}
\sigma^1_{1} \sigma^1_{2} ...   \sigma^1_{n-1} \sigma^1 $$
\be \label{S12n}  
  S^{12}_n = 
S^{12}(2\ell1_n-2\ell^1_1) S^1_{n-1}(2\ell^1_n-2\ell^1_2)
S^1_{n-2}(2\ell^1_n-2\ell^1_3) ...S^1_{2}(2\ell^1_n-2\ell^1_{n-1}) $$ $$
S^1_1(2\ell^1_n-2\ell^2_n) S^1_2(2\ell^1_{n-1}-2\ell^2_n) ...
S^1_{n-1}(2\ell^1_{2}-2\ell^2_n) S^{12}(2\ell^1_1-2\ell^2_n) 
\ee
It is a 1 to 1   beta sequence,
$$  S^1_n \cdot 1 = \pi^1_n(\mathfrak l) 1, $$
$$ \pi^1_n(\mathfrak l) =
\prod_1^{n-1} B(2\ell^2_n-2\ell^1_{k+1}, 2\ell^1_k- 2\ell^2_n +1)
\prod_1^{n-2} B(2\ell^1_{k+1}-2\ell^1_n, 2\ell^1_k-21\ell^1_{k-1}+1). $$

General beta sequences for $g\ell(n)$ type Yangian presentations will be
derived  in section 9. In the next section we consider the details of the
$g\ell(2)$ case and show how the permutation coefficients distinguish 
the representation types.


\section{On the Yangian representations of type $g\ell(2)$ }
\setcounter{equation}{0}

\subsection{Explicit  $L$ operators}

Let us write the Jordan-Schwinger $L^-$ form of the $g\ell(2)$ generators as the
matrix
\be \label{L0gl2} (L_{ab}) = -
\begin{pmatrix}
\dd_1 x_1& \dd_1 x_2 \\
\dd_2 x_1 & \dd_d x_2
\end{pmatrix}.
\ee 
We reduce the action  to 
 functions of the form
$(x_2)^{2\ell_2} \ \varphi(x) $,
$x= \frac{x_1}{x_2} $. 
In particular, this leads to
$$ \dd^1_2 x^1_2 (x_2^1)^{2\ell_2} \ \varphi(x^1) 
=  (x_2^1)^{2\ell_2}  (2\ell_2 + 1 + x^1_2 \dd^1_2 x^1_2 ) \varphi
=  (x_2^1)^{2\ell_2}  (2\ell_2 + 1 -  x^1_1 \dd^1_1 x^1_1 ) \varphi
= $$ $$
(x_2^1)^{2\ell_2}  (2\ell_2 +2 - \dd^{1 \prime }_1 x^{1 \prime}_1 ) \varphi.
$$
We obtain $L(u^-,u), u^-=u-2\ell$  (\ref{Lu-u})
 in matrix form in terms of the normal
coordinates $x =\frac{x_1}{x_2}, \dd = x_2 \dd_1 $ 
\be \label{Lugl2} 
(L_{ab}(u^-,u) ) =
 \begin{pmatrix}
u- 1 -x \dd& -\dd \\
-x(u-u^--1- x \dd) & u^--1 + x \dd 
\end{pmatrix}.
\ee 

 The quantum determinant is
$$ qdet L(u^-,u) = (u^--2) (u-1), $$.


In the case of rank 1 all highest weight representations can be formulated in the JS
presentation. The Biedenharn construction is actually not  needed here. 
It is done in one step and this step supplies merely the second representation
parameter $2\ell_1$.  The canonical pairs of the second factor $L^2$ are
completely removed in the procedure.     

We start from the two-factor monodromy in the form (\ref{L1L2}) and impose the
constraints  (\ref{C2}), (\ref{gauge}).   
We start with the second factor $L^2(0)$ obtained from (\ref{L0gl2}) by
indexing the canonical pairs with superscript 2.
The constraint $C^2$ 
allows to eliminate the pair $x^2_2, \dd^2_2 $,

$$ \dd^2_2 = - \frac{1}{x^1_2 } x^1_1 \dd^2_1. $$
The constrained factor $L^2$ has the form
$$ - L^2(0)|_{C^2} =
\begin{pmatrix}
\dd^2_1 x^2_1 & 0  \\
  - \frac{1}{x^1_2 } x^1_1 \dd^2_1 x^2_1 & 0   \\   
\end{pmatrix}.
$$
$L^1(0)$ is obtained from (\ref{L0gl2}) by
indexing the canonical pairs with superscript 1.

We check that $ L^1(0) L^2(0)|_{C^2} = 0 $ holds after this step. 
 
We introduce the  projection on homogeneous functions.
We are careful to keep the ordering  prescription
adopted for the  constraint.

In the case $L^1$ we reduce like as above  to functions of the form 
$(x_2^1)^{2\ell_2} \ \varphi(x^{1 \prime}_1) $,
$x^{1\prime }_1= \frac{x^1_1}{x^1_2} $. 
$$ L^1((u-2\ell_2,u) =
\begin{pmatrix}
u-\dd^1_1 x^1_1 & -\dd^1_1 x^1_2  \\
 -(2\ell_2 + 2 - \dd^1_1 x^1_1) \frac{ x^1_1}{ x^1_2} 
& u-  (2\ell_2 + 2 - \dd^1_1 x^1_1)   \\   
\end{pmatrix} =
$$ 
\be \label{L10}
\begin{pmatrix}
u-\dd^{1 \prime}_1 x^{1 \prime}_1 & -\dd^{1 \prime}_1   \\
 -(2\ell_2 + 2 - \dd^{1 \prime}_1 x^{1 \prime}_1)  x^{1 \prime}_1 
& u- (2\ell_2 + 2 - \dd^{1 \prime}_1 x^{1 \prime}_1)   \\   
\end{pmatrix}. 
\ee
In the case $L^2$
 the reduction to functions of the form $ 
(x_1^2)^{2\ell_1} \varphi $ 
leads  to
$$ \dd^2_1 x^2_1 (x_1^2)^{2\ell_1} \ \varphi
=  (x_1^2)^{2\ell_1}  (2\ell_1 + 1 + x^2_1 \dd^2_1  ) \varphi
=  (x_1^2)^{2\ell_1}  (2\ell_1 + 1 -  x^2_2 \dd^2_2  ) \varphi.
$$
The constraint $C^2$ implies that $\varphi$ is constant and the
last term in the bracket drops out. With this line of argument 
the constraint and the  projection commute
and result in 
\be \label{L20}
 L^2(u-2\ell_1, u)|_{C^2} =
\begin{pmatrix}
u-2\ell_1 -1 & 0   \\
 + x^{1 \prime}_1 (2\ell_1+1) 
& u   \\   
\end{pmatrix}.
\ee
We check the relation 
\be \label{L1L20}
 L^1(u-2\ell_2, u) L^2(u-2\ell_1,u)|_{C^2} = u \mathcal{L}^{1+2} (u). \ee
 We have obtained
in terms of normal coordinates (primes omitted)
\be \label{L1+2mat}
\mathcal{L}^{1+2}(u; 2\ell_1, 2\ell_2) = 
\begin{pmatrix}
u-2-2\ell_1- x \dd & -\dd \\
-x (2\ell_2 -2\ell_1-1- x \dd) & u-2\ell_2-1+ x \dd 
\end{pmatrix}.
\ee
We calculate the quantum determinant 
$$ qdet L^{1+2}(u;2\ell_2,2\ell_1) = $$ $$
(u-2\ell_1 -2 - x^{1 }_1 \dd^{1 }_1) 
(u-2\ell_2-1-1+ x^{1 }_1 \dd^{1 }_1 ) -
x^{1 }_1(2\ell_2-2\ell_1 -1 - x^{1 }_1\dd^{1 }_1)
\dd^{1 }_1  =
$$ \be \label{qdet01}
(u-2\ell_1 -2) (u-2\ell_2 -2).
\ee

\subsection{Permutations }

We consider the second order evaluation, $N=2$, where
the generators $T_{ab}(u;2\frak l)$  (\ref{monogen}) are obtained as the matrix  
product of two factors $\mathcal {L}^I(u,2\underline \ell^I)$, 
(\ref{Lgen}), explicitly written in (\ref{L1+2mat}),  
\be \label{mono2}
T(u;2\frak l) = \mathcal {L}^1 (u,2\underline \ell^1)\mathcal {L}^2(u,2\underline
\ell^2) . \ee
The parameter array is
$$ 2\frak l = 2\underline \ell^1; 2\underline \ell^2 =
2\ell^1_2, 2\ell^1_1; 2\ell^2_2, 2\ell^2_1. $$

Consider the permutations   
\be \label{sigma12}
\sigma^{12}_1: 2\ell^1_1 \leftrightarrow 2\ell^2_1,  \ \ \
\sigma^{12}_2: 2\ell^1_2 \leftrightarrow 2\ell^2_2 \ee
permuting $2 \ell^I_1$ and $2\ell^I_2$, respectively. The
 permutation $\sigma^{12} = \sigma^{12}_1\sigma^{12}_2 $
 exchanges  
$2 \underline \ell^1$ with  $ 2\underline \ell^2$.

They can be obtained as products of the adjacent permutations $\sigma^1,
\sigma^I_1,  I=1,2 $ (\ref{sigmaI})
$$ \sigma^{12}_1= \sigma^1\sigma^2_1 \sigma^1, \ \ 
\sigma^{12}_2= \sigma^1\sigma^1_1 \sigma^1. $$
The corresponding sequences of permutation operators acting on the monodromy
are
$$  S^{12}_1= S^1(2\ell^2_2-2\ell^2_1)  S^2_1(2\ell^2_1- 2\ell^1_1)
S^1(2\ell^1_1-2\ell^2_2), $$
\be \label{S122}  
S^{12}_2= S^1(2\ell^1_2-2\ell^1_1) S^1_1(2\ell^2_2-2\ell^1_2) 
S^1(2\ell^1_1-2\ell^2_2).
\ee
We calculate the action on $1$. We write the details for the latter case.
$$ S^1(2\ell^1_1-2\ell^2_2) \cdot 1 = (x^1-x^2)^{2\ell^1_1-2\ell^2_2}, $$
$$ S^1_1(2\ell^2_2-2\ell^1_2)S^1(2\ell^1_1-2\ell^2_2) \cdot  1 =
\int dc c^{2\ell^2_2-2\ell^1_2-1} (x^1-x^2-c)^{2\ell^1_1-2\ell^2_2}.
 $$ 
The operator is specified by defining the integration path
as the Pochhammer contour  encompassing the
two branch points $c=0$ and $c= x^1-x^2$. We use the integral representation of the Euler
beta function  and obtain 
$$
B(2\ell^2_2-2\ell^1_2, 2\ell^1_1-2\ell^2_2+1) \ 
(x^1-x^2)^{2\ell^1_1-2\ell^1_2}. $$
The action of  $ S^1(2\ell^1_2-2\ell^1_1)$   amounts in the multiplication by
$(x^1-x^2)^{2\ell^1_2-2\ell^1_1}$.

We obtain
$$  S^{12}_i (2\frak l) \cdot 1= \pi^{12}_i (2\frak l) \cdot 1, $$
\be \label{pii}
 \pi_i (2\frak l) = \pi_i (2 \underline \ell^1, 2 \underline \ell^2) = 
B(2\ell^2_i-2\ell^1_i, 2\ell^1_1-2\ell^2_2+1). \ee

The operator representing $\sigma^{12} $ in action on the monodromy
$T(u; 2\frak l)$ is
\be \label{S^12}
  S^{12} (2\frak l) =  S^{12}_2 (\sigma^{12}_2 2\frak l)
 S^{12}_1 (2\frak l) =  S^{12}_1 (\sigma^{12}_1 2\frak l)
 S^{12}_2 (2\frak l) = \ee $$
S^1(2\ell^2_2-2\ell^2_1)  S^2_1(2\ell^2_1- 2\ell^1_1)
S^1_1(2\ell^2_2-2\ell^1_2) S^1(2\ell^1_1-2\ell^2_2).
$$
Its action on $1$ is
$$  S^{12} (2\frak l) \cdot 1= \pi (2\frak l) \cdot 1, $$
\be \label{pi}
 \pi (2\frak l) = \pi (2 \underline \ell^1, 2 \underline \ell^2) =
\Gamma(2\ell^2_2 -2\ell^1_2) \Gamma(2\ell^2_1 -2\ell^1_1)
\frac{\Gamma(2\ell^1_1-2\ell^2_2+1)}{ \Gamma(2\ell^2_1-2\ell^1_2+1)}.
\ee  

In analogous steps   we can calculate the action on
monomials $(x^1)^{m_1} (x^2)^{m_2}$, 
$$
  S^{12}_i (2\frak l) \cdot (x^1)^{m_1} (x^2)^{m_2} =  
(x^1)^{m_1} (x^2)^{m_2}  \ \ \ \Gamma( 2\ell^1_1-2\ell^2_2+1) \cdot $$ 
\be \label{S12im}
\sum_k^{m^{i\p}}  
\begin{pmatrix}
m^{i\p} \\ k 
\end{pmatrix}
(\frac {x^i }{x^{i\p} } -1)^k  
\frac{\Gamma (2\ell^2_i-2\ell^1_i+k)}
{\Gamma(2\ell^{i\p}_1-2\ell^{i\p}_2+1+k) }.
\ee
Here we use the abbreviation $i\p = 3-i$. 

We notice that the point permutation
$$\mathcal P^{12} (x^1, \dd^1) = (x^2, \dd^2) \mathcal P^{12} $$
together with $ S^{12} $ results in the monodromy with the reversed order
of the factors,
$$  S^{12} (2\frak l) T(u, 2\frak l)  ( S^{12} (2\frak l) )^{-1}
= T(u, \sigma^{12} 2 \frak l), $$ $$
\mathcal P^{12} T(u, \sigma^{12} 2 \frak l)  \mathcal P^{12} =
\mathcal L^2(u, 2\underline \ell^2) \mathcal L^1(u, 2\underline \ell^1). $$
The latter algebra is equivalent to the algebra with the generators
$T(u,2\frak l)$. 
This known fact is confirmed here by observing  that the representation 
vectors of generated by $T(u, \sigma^{12} 2 \frak l)$ on $ 1$ differ from the corresponding
ones of 
$T(u,  2 \frak l) $ on $ 1$  by the exchange of $ x^1, x^2 $.

\subsection{Representations of the second order evaluation}

The monodromy matrix of second order $N=2$ (\ref{mono2}) decomposes in $u$ as
$$ T_{ab}(u; 2\ell^1_2, 2\ell^1_1, 2\ell^2_2, 2\ell^2_1) = 
u^2 \delta_{ab} + u T_{ab}^{[1]} + T_{ab}^{[2]}  $$
The generators can be calculated explicitly by (\ref{L1+2mat}).

$$T_{12}^{[1]} = - \dd^1 - \dd^2, $$
$$ T_{12}^{[2]} = (2+2\ell^1_1 + x^1 \dd^1) \dd^2 + \dd^1 (1+2\ell^2_2
-x^2\dd^2) = $$
$$(\dd^1+\dd^2) 
\half [3+2\ell^1_1 +2\ell^2_2 + x^1\dd^1-x^2 \dd^2]  -
(\dd^1-\dd^2) 
\half [1+2\ell^1_1-2\ell^2_2 +x^1\dd^1 +x^2 \dd^2 ], 
$$ 

$$T_{21}^{[1]} = -x^1(2\ell^1_2-2\ell^1_1-1 -x^1\dd^1) -
x^2(2\ell^2_2-2\ell^2_1-1 -x^2\dd^2),
$$ 
$$
T_{21}^{[2]} = x^1 (2\ell^1_2-2\ell^1_1-1 -x^1\dd^1) (2+2\ell^2_1+x^2\dd^2) 
+ x^2(2\ell^2_2-2\ell^2_1-1 -x^2\dd^2) (2\ell^1_2 +1 - x^1\dd^1) 
=$$ 
$$ \half [x^1(2L^1-x^1\dd^1) -x^2(2L^2-x^2\dd^2)] (1+2\ell^2_1-2\ell^1_2+
x^1\dd^1+x^2\dd^2 ) + $$ $$
\half [x^1(2\ell^1_2-2\ell^1_1-1-x^1\dd^1) +
 x^2(2\ell^2_2-2\ell^2_1-1-x^2\dd^2)] (3+2\ell^2_1 +2\ell^1_2 +x^2\dd^2
-x^1\dd^1 ). $$

We look for a highest weight vector. The condition $ T_{12}^{[1]} \psi = 0$ leads to
the ansatz $ \psi_0 = \psi (x^1-x^2)$. The remaining condition
$$ T_{12}^{[2]} \psi (x^1-x^2) = [(2\ell^2_2 -2\ell^1_1 -1) - (x^1-x^2)\dd^1 ]
\dd^1 \psi (x^1-x^2) = 0 $$
is solved besides of $ \psi = 1$ by
\be \label{M12}
 \psi_0 (x^1-x^2) =  (x^1-x^2)^{2M^{12}+1}, \ \ \ 
2M^{12} = {2\ell^2_2 -2\ell^1_1 -1}. \ee
The representation generated by $T_{ab}(u)$ on $\psi_0$ is spanned by
monomials in $x^1, x^2$ multiplied by $\psi_0$ with the weight functions
$$\tilde \lambda_1(u) = (u-2-2\ell^2_2) (u-2-2\ell^2_1), \ \  
\tilde \lambda_2(u) = (u-2-2\ell^1_2) (u-2-2\ell^1_1). $$

We find a lowest weight vector $\psi_-$  obeying $ T_{21}(u) \psi_- = 0 $, 
$$  \psi_- = (x^1)^{ 2L^1} (x^2)^{ 2L^2}, $$
\be \label{L1,L2}
 2L^1 = 2\ell^1_2-2\ell^1_1-1 =2L^1 , 2L^2= 2\ell^2_2-2\ell^2_1-1. \ee

If both $2L^1, 2L^2 $ are non-negative
integers, then  $\psi_-$ lies in the space spanned by the monomials in
$x^1, x^2 $ and     
 the rising operation starting 
from the highest weight vector $1$ terminates at $\psi_-$.   

In the case of generic parameters the results of the action of $T_{21}^{[1]} $ and $T_{21}^{[2]} $
on a monomial $ (x^1)^{ m_1} (x^2)^{ m_2} $ are linear independent. 
We look for the  particular parameter values, where  the rising operation degenerates at
some values of the exponents $m_1,  m_2$, in the way that
the result of the action of $T_{21}^{[1]} $ and $T_{21}^{[2]} $ 
are linear dependent,
\be \label{AB} 
A T_{21}^{[1]} (x^1)^{ m_1} (x^2)^{ m_2} + B T_{21}^{[2]} (x^1)^{ m_1} (x^2)^{ m_2} =
0 ,\ee 
$$ T_{21}^{[1]} (x^1)^{ m_1} (x^2)^{ m_2} = 
(x^1)^{ m_1} (x^2)^{ m_2} [ -x^1 ((2\ell^1_2-2\ell^1_1-1 -m_1) - x^2 (2\ell^2_2-2\ell^2_1-1
-m_2) ], $$ 
$$ T_{21}^{[2]} (x^1)^{ m_1} (x^2)^{ m_2} = (x^1)^{ m_1} (x^2)^{ m_2}
[x^1 ((2\ell^1_2-2\ell^1_1-1 -m_1)(2+2\ell^2_1+m^2) + $$ $$
x^2 (2\ell^2_2-2\ell^2_1-1
-m_2) (2\ell^1_2+1-m_1)].
$$
We decompose the condition with respect to the independent variables $x^1,
x^2$. The system of homogeneous equations has a non-trivial solution in
$A, B$ if the determinant vanishes. We obtain the degeneracy condition
$$(2\ell^1_2-2\ell^1_1-1 -m_1) (2\ell^2_2-2\ell^2_1-1-m_2) 
(2\ell^1_2-2\ell^2_1-1-m_1-m_2) = 0. $$

The degeneracy may appear, if there is a non-negative integer value of
$2L^1, 2L^2$ or of
\be \label{M21}
 2 M^{21} = 2\ell^1_2-2\ell^2_1-1, \ee
correspondingly at the monomials
$$ (x^1)^{2L^1} (x^2)^{m_2}, \ \ (x^1)^{m_1} (x^2)^{2L^2}, \ \ (x^1)^{m} (x^2)^{2M^{21}-m} , m <
2M^{21}+1.
$$
Comparing the forms of the equivalent algebras of $T(u, 2\frak l)$ and
$T(u, \sigma^{12} 2\frak l)$, we notice that
$\sigma^{12} 2L^1 = 2L^2, \sigma^{12} 2M^{12} = 2M^{21} $ and thus 
the role of $\psi_0 $ and 
$$ \psi_1 = (x^1-x^2)^{2M^{21} +1}.$$
 is interchanged.  

We show that the permutation operations considered in subsection 8.2 
allow to derive the above results without using the explicit form of the
generators. 

The existence of a highest weight vector besides of the constant $1$ in the
space of functions of $x^1, x^2$ can be derived starting from
$T(u, \sigma^1 2 \frak l), \sigma^1 2 \frak l = 2\ell^1_2,2\ell^2_2;2\ell^1_1,
2\ell^2_1 $,   
$$ T_{12}(u, \sigma^1 2 \frak l) \cdot 1 = 0, \ \ \
T_{aa}(u, \sigma^1 2 \frak l) \cdot 1 = \lambda_a(u, \sigma^1 2 \frak l)
\cdot 1. $$
We act by $S^1(2\ell^2_2-2\ell^1_1$ to obtain
$$ T_{12}(u,  2 \frak l) \cdot \psi_0(x^1,x^2)  = 0, \ \ \
T_{aa}(u,  2 \frak l) \cdot \psi_0(x^1,x^2) = \lambda_a(u, \sigma^1 2 \frak l)
\psi_0(x^1,x^2) $$
$$ \psi_0(x^1,x^2) = S^1(2\ell^2_2-2\ell^1_1) \cdot 1 =
(x^1-x^2)^{2\ell^2_2-2\ell^1_1} $$
This confirms that $\psi_0$ is a highest weight vector of the algebra $T(u,  2 \frak
l)$. 
 
The representation generated by $T(u, \sigma^1 2 \frak l)$
on $1$ is mapped to the representation generated by 
 $T(u,  2 \frak l) $ on $\psi_0(x^1,x^2)$.
In particular we see that the roles played by $2L^1, 2L^2, 2M^{12}$ 
in the form of the algebra $T(u, \sigma^1 2 \frak l)$ are
 $$\sigma^1 2L^1 = 2\ell^1_2-2\ell^2_2-1, \ \ \sigma^1 2L^2 =
2\ell^1_1-2\ell^2_1-1, \sigma^1 2M^{12} = 2\ell^1_1 -2\ell^2_2-2,  $$
$$ \sigma^1 2 L^1 = -2M^{12} + 2L^1-1, \sigma^1 2 L^2 = -2M^{12} + 2L^2-1, $$
\be \label{sigma1M}
\sigma^1 2M^{12} +1 = -( 2M^{12} +1), \sigma^1 2M^{21} = 2M^{21}. \ee
 
In analogy we confirm the role of $\psi_1$ as a highest weight vector of
$T(u,\sigma^{12} 2\frak l)$. Starting from $\tilde \sigma = \sigma^1
\sigma^{12} = 2\ell^2,2\ell^1_2; 2\ell^2_1,2\ell^1_1 $, 
$T(u, \tilde \sigma 2\frak l) \cdot 1$ and acting by 
$S^1(2\ell^1_2-2\ell^2_1) $ we obtain
$$ T_{12}(u, \sigma^{12}  2 \frak l) \cdot \psi_1(x^1,x^2)  = 0, \ \ \
T_{aa}(u, \sigma^{12} 2 \frak l) \cdot \psi_1 (x^1,x^2) = \lambda_a(u,
\tilde \sigma 2 \frak l)
\psi_1(x^1,x^2), $$
$$ \psi_1(x^1,x^2) = S^1(2\ell^2_1-2\ell^1_2) \cdot 1 =
(x^1-x^2)^{2\ell^1_2-2\ell^2_1}. $$  
The roles of $2L^1, 2L^2, 2M^{12}, 2 M^{21} $ are played by
$$ \tilde \sigma 2L^1= 2\ell^2_2-2\ell^1_2-1, \ \ 
\tilde \sigma 2L^2= 2\ell^2_1-2\ell^1_1-1,  \ \ 
\tilde \sigma 2M^{12} = 2\ell^2_1-2\ell^1_2-1,  $$
$$ \tilde \sigma 2L^1= 2M^{12} -2L^1-1, \tilde \sigma 2L^2= 2M^{12}
-2L^2-1,
$$ \be \label{sigmatM}
\tilde \sigma 2M^{21}  = 2M^{12} , \tilde \sigma 2M^{12} +1 = -(2M^{21} +1). \ee

Consider the finite-dimensional representation generated by  $T(u,  2 \frak l) $ 
on $1$ for non-negative integer  
$2 L^I$. It is spanned by $(x^1)^{m_1}  (x^2)^{m_2} , 
m_1 \le 2 L^1, m_2 \le 2 L^2 $.

We have a second highest weight vector in this representation space 
$$ \psi_0 =( x^1-x^2)^{2M^{12}+1}, 2M^{12}=2\ell^2_2 - 2\ell^1_1 -1 $$
if $ 2M^{12}  $ is non-negative integer and if
$$ 2 L^I > 2M^{12} \ge 0 , I=1,2. $$
This means
\be \label{l22<0}
 2\ell^2_2 - 2\ell^1_2 < 0, \ \ 2\ell^2_1 - 2\ell^1_1 < 0. 
\ee
The parameter combinations introduced above obey
\be \label{LL=MM}
 2L^1+2L^2= 2M^{12} +2M^{21}. \ee
If $2 L^1, 2 L^2, 2M^{12} $ are integers then $2M^{21}$ is integer. 
We have seen that the
 rising operation by the generators $T_{21}(u) $ 
on  monomials $(x^1)^{m_1}  (x^2)^{m_2} $ degenerates 
in the case of non-negative integer values of $2M^{21}$.
We consider the condition that for integer $2 L^1, 2 L^2, 2M^{12} $
the latter barrier lies in the finite-dimensional
subspace, 
$$ 2L^1 \ge 2M^{21}, 2 L^2 \ge 2M^{21}. $$
This means
\be \label{l22>0}
 2\ell^2_2 - 2\ell^1_2 \ge  0, \ \ 2\ell^2_1 - 2\ell^1_1 \ge  0. \ee 
This is the opposite of both conditions (\ref{l22<0}).

The finite dimensional representation  is irreducible and  spanned by  all the
monomials with powers bounded by $2L^1, 2L^2 $ if 

$$ 2\ell^2_2 - 2\ell^1_2 < 0, 2\ell^2_1 - 2\ell^1_1 \ge 0 $$
 or if 
$$ 2\ell^2_2 - 2\ell^1_2 \ge 0, 2\ell^2_1 - 2\ell^1_1 < 0. $$

Thus the values of the parameter combinations (\ref{M12}), (\ref{L1,L2}),
(\ref{M21}),
$2L^1, 2L^2, 2M^{12}, 2M^{21} $ and the ordering relations among them
defines the type of representation. 

Let us consider the action of the parameter permutations 
$\sigma^{12}_1, 
\sigma^{12}_2$ (\ref{sigma12})
on the considered combinations. 
$$ \sigma^{12}_1: M^{12} \leftrightarrow 2 L^2,  M^{21} \leftrightarrow 2
L^1, 
$$
$$ \sigma^{12}_2: M^{12} \leftrightarrow 2 L^1,  M^{21} \leftrightarrow 2
L^2. $$
The ordering configurations 
\be \label{M12LLM}
 2M^{12} \le 2L^1, 2L^2 \le 2M^{21} \ee are mapped 
by $\sigma^{12}_1$ to 
\be \label{L2MML}
 2L^2 \le 2M^{12}, 2M^{21} \le 2L^1, \ee 
 by $\sigma^{12}_2$ to 
\be \label{L1MML}
 2L^1 \le 2M^{12}, 2M^{21} \le 2L^2, \ee 
and by $\sigma^{12}$ to
\be \label{M21LLM} 
 2M^{21} \le 2L^1, 2L^2 \le 2M^{12}. \ee

The equivalence of the algebras $T(u,2\frak l)$ and $T(u,\sigma^{12}\frak l)
$ allows the restriction to two configurations, (\ref{L1MML}) and
(\ref{M12LLM}). 
In the case
\be \label{L1<MM<L}
 2L^1<2M^{12}, 2M^{21} < 2L^2 \ee
the representation generated by $T(u,2\frak l)$ on $1$ is finite dimensional
irreducible and the representation generated by $T(u,2\frak l)$ on $\psi_0$  
is infinite dimensional with degeneracy of the rising operation. 

In the case
\be \label{M12<LL<M} 
2M^{12} < 2L^1, 2L^2 < 2M^{21} \ee
the representation generated by $T(u,2\frak l)$ on $1$ is finite dimensional
but not irreducible, because $\psi_0$ lies inside the finite-dimensional subspace.
The representation generated by $T(u,2\frak l)$ on $\psi_0$ is finite dimensional
irreducible.

In the limiting case 
\be \label{L1=M12}
 2L^1= 2M^{12}, 2L^2 = 2M^{21}, 2L^2 < 2M^{12}  \ee
the representation generated by $T(u,2\frak l)$ on $1$ is finite dimensional
irreducible and the representation generated by $T(u,2\frak l)$ on $\psi_0$
is infinite dimensional with degeneracy in the rising operations increasing
$m_2$. 

In the limiting case 
\be \label{L1=L2=M12}
 2L^1= 2L^2 = 2M^{12}= 2M^{21}  \ee
the representation generated by $T(u,2\frak l)$ on $1$ is finite dimensional
irreducible and the representation generated by $T(u,2\frak l)$ on $\psi_0$ 
is irreducible infinite dimensional.

\subsection{Representation type by permutation coefficients}

Let us write the permutation coefficients (\ref{pii}) in terms of the parameter
combinations $2L^1, 2L^2 $, $2M^{12}, 2M^{21} $ (\ref{M12}),
(\ref{L1,L2}), (\ref{M21}). We have discussed how  their special values
and ordering relations characterize the  representation types. We consider
now how this appears in zeros and poles of the permutation coefficients.

\be \label{piiM}
\pi_1 (2\frak l) = B(2M^{12}-2L^1, -2M^{12}), \ \  
\pi_2 (2\frak l) = B(2M^{12}-2L^2, -2M^{12}). \ee
They vanish correspondingly for non-negative integer values of
$2L^i$. 

We notice by (\ref{S12im} ) that in the case of non-negative integer $2L^i$ 
the action  of $ S^{12}_{i\p}, (i\p+i=3 ) $  on monomials vanishes if
$ m_i < 2 L^i+1 $. If both $2L^1, 2L^2$ are non-negative integers then
a projection on the subspace spanned by the monomials with 
$ m_1 <2L^1+1, m_2 <2L^2+1 $ can be constructed.
  Notice that in $\pi(2 \frak l)$ the $(\Gamma (-2L^i))^{-1}$
in the second factor is canceled by  its inverse appearing in the first factor.
Indeed, 
$$ \pi_2 (\sigma^{12}_1 2\frak l) = B (2L^2-2M^{12}, -2L^2 ),  
 $$
\be \label{piM}
\pi(2 \frak l) = 
 \pi(2\underline \ell^1, 2\underline \ell^2) = 
\Gamma(2M^{12} - 2L^1) \Gamma(2M^{12} -2L^2) \frac{\Gamma(-2M^{12}
)}{\Gamma(-2M^{21}) }. \ee
  
The representation generated by $T(u; \sigma^{12}_1 2\frak l)$ on $ 1 $ has weight
functions coinciding with the ones of the representation generated by $T(u;  2\frak
l)$ on $ 1$. For generic parameters both representations are irreducible and
equivalent. If $2L^2$ is non-negative integer then the latter is not irreducible,
the rising operation $T_{21}(u)$ does  generate monomials $(x^1)^{m_1}
(x^2)^{m_2}   $ with 
$m_2 <2L^2+1$ only. The permutation operator  $ S^{12}_1 (2\frak l) $
mapping to the first one 
annihilates these monomials (\ref{S12im} ). In the image of $\sigma^{12}_1$ the value of $2L^2$ becomes
 the value of $\sigma ^{12}_1 2L^2 = 2M^{12} $ and  here  the representation is characterized by the
appearance of a second highest weight vector. This is indicated by the
vanishing of the permutation factor $ \pi_2 (\sigma^{12}_1 2\frak l)$
(\ref{piM}).

The first arguments in $\pi_i (2\frak l)$ or $ \pi_i (\sigma^{12}_{i\p} 2\frak l) $   
are related to the degeneracy of the rising operation of $T(u,2\frak l)$ on
$\psi_0$ or $T(u, \sigma^{12} 2 \frak l)$ on $\psi_1$, as obtained in
(\ref{sigma1M}), (\ref{sigmatM}).
If $2M^{12}-2L^I $ is
integer  the powers of $x^I$ are limited
for negative values in the first case and for positive values of these arguments in the second
case, correspondingly. No limits are encountered in the cases
of $2M^{12}-2L^I =0 $. 
If both $2L^2$ and $2M^{12}$ are non-negative integers the subspace spanned
by the monomials $m_2 <2L^2+1$ carries the irreducible representation of 
  $T(u;  2\frak l) $ on $1$ in the case $2 L^2 < 2M^{12}$. 
In the opposite
case the representation on
this subspace is not irreducible, containing a second highest weight vector.        

Now we consider parameter arrays with a shift $u$ and show that the
asymptotics of permutation coefficients at $u\to 0$ distinguishes the
representation types.

If $2\ell^2_a = 2\ell^1_a+u$ then the monodromy of second order $T(v;2\frak l)$ 
(\ref{mono2}) is composed of equivalent   $g\ell(2)$
generators. Let us consider the $u$ dependence of the resulting permutation
factors and compare it to the features of the representations.
$$ \pi_1(2\underline \ell^1, 2\underline \ell^1 +u) = 
\pi_2(2\underline \ell^1, 2\underline \ell^1 +u) = B(u,-2L^1-u),
$$
\be \label{pil1}
\pi(2\underline \ell^1, 2\underline \ell^1 +u) = (\Gamma(u))^2 
\frac{\Gamma(2L^1-u)}{ \Gamma(2L^1+u)}. \ee  
It is well known that  in the case of generic parameters
the representation of $\mathcal L(u, 2\ell^1_2,2\ell^1_1) $ with the highest
weight
vector $1$ is irreducible  and  spanned by all monomials $x^m$. If  
$2L = 2\ell_2-2\ell_1-1$ is non-negative integer the irreducible
representation is spanned by the monomials with $m<2L +1 $.
The two types of representations are distinguished in the $u \to 0$
asymptotics of the permutation coefficient.
In the generic case 
$$ \pi(2\underline \ell^1, 2\underline \ell^1 +u) \sim u^{-2}, $$
but in the case of non-negative integer $2L^1$ the leading term has the
opposite sign,
$$\pi(2\underline \ell^1, 2\underline \ell^1 +u) \sim - u^{-2}. $$   

Let us write explicitly also the permutation factors resulting by the
substitution $2  \ell_a^2 \to 2 \ell^2_a +u$ and compare the behavior at
$u\to 0$ to the 
cases of the representations of the second order evaluation. 
$$ \pi_1(2\underline \ell^1, 2\underline \ell^2 +u) =
 B(2M^{12}-2L^2+u, -2M^{12}-u),  \ \ 
\pi_2(2\underline \ell^1, 2\underline \ell^1 +u) = 
B(2M^{12}-2L^1+u,-2M^{12}-u), $$
\be \label{pil1l2} 
\pi(2\underline \ell^1, 2\underline \ell^2 +u) = 
\Gamma(2M^{12}-2L^1 +u) \Gamma(2M^{12}-2L^2 +u) \frac{\Gamma
(-2M^{12}-u)}{\Gamma(-2M^{21}+u) } = \pi^{12}(u) \ee
Also in these expressions the asymptotics in $u$ changes 
at non-negative integer $2 M^{12}, 2L^1, 2L^2$ compared to the 
one at generic values. 

Consider the coefficient $\pi(2\underline \ell^1, 2\underline \ell^2 +u)$
(\ref{pil1l2}) in the case of non-negative integer 
$2L^1, 2L^2$. For generic $2M^{12}$ the expansion in  $u$ starts with 
a non-vanishing constant. 
If $2M^{12}$ is  integer too
then we have first the case  of a negative integer
$2M^{12}$ in the configuration (\ref{M12LLM}), where the asymptotics is  constant times $ (-1)^{2M^{12}+1}
u^{-3} $. 

In the cases of all $2L^1, 2L^2, 2M^{12},2M^{21} $ 
non-negative integer and not coinciding we have the asymptotics
correspondingly in the configurations 
\begin{enumerate}
\item
 (\ref{M12LLM}) $\pi^{12}(u) \sim - u^{-2}$, 
\item
(\ref{L1MML})  $ \pi^{12}(u) \sim (-1)^{2M^{12}+2L^2+1}
u^{-1}, $
\item
 (\ref{L2MML})  $ \pi^{12}(u) \sim (-1)^{2M^{12}+2L^1+1}
u^{-1}, $
\item
 (\ref{M21LLM})  $ \pi^{12}(u) \sim (-1)^{2L^1+2L^2+1}. $
\end{enumerate}

The asymptotics of the first case applies also for $2M^{12}= 2L^1\le 2L^2$
and $2M^{12} =2L^2 \le 2L^1$. The second and the third cases apply
correspondingly, if $2L^1 < 2M^{12}=2M^{21} $ or $2L^2 < 2M^{12}=2M^{21} $.
The fourth case applies also for $2M^{21} < 2L^1=2L^2$.

We consider the the monodromy 
$$ T_{4 }(v;2 \frak l^{12}, \frak l^{34}) = \mathcal L^1(v;2\ell^1_2,2\ell^1_1)
\mathcal L^2(v;2\ell^2_2,2\ell^2_1)\mathcal L^3(v;2\ell^I_3,2\ell^3_1)\mathcal
L^4(v;2\ell^4_2,2\ell^4_1)
$$
of the 
$N=4$ evaluation of $\mathcal{Y}(g\ell(2))$
characterized by the parameter array
$$ 2 \frak l^{12}, \frak l^{34}= 2 \underline \ell^1,   2 \underline \ell^2, 2 \underline \ell^3, 2 \underline
\ell^4. $$
The monodromy is constructed as the matrix product of the $L$ matrices $\mathcal
L^I(v;2\ell^I_2,2\ell^I_1)$ obtained from (\ref{L1+2mat}) by substituting the
canonical pair $x, \dd$ by $ x^I, \dd^I$ and the parameters $2\ell_2, 2\ell_1
$ by $2\ell^I_2,  2\ell^I_1$ .

In order to study the reducibility in the $N=2$ evaluation we construct the
permutation $S^{1234} $ to the algebra with the parameters
$$ 2 \frak l^{34}, 2 \frak l^{12} = 2 \underline \ell^3,   2 \underline \ell^4, 2 \underline \ell^1, 2 \underline
\ell^2. $$
This parameter permutation $\sigma^{1234}$ can be constructed from the permutations
$\sigma^{I,I+1}$
transforming the sequence $2 \underline \ell^I, 2\underline \ell^{I+1} $ to
$2\underline \ell^{I+1}, \underline \ell^I$ in analogy to the particular
case $ \sigma^{12}, I=1 $ above (\ref{sigma12}),
$$ \sigma^{1234} =  \sigma^{23} \sigma^{34} \sigma^{12} \sigma^{23}. $$
The corresponding operator acting on the monodromy is the product of 
$ S^{I,I+1}$ obtained from $ S^{12}$ (\ref{S12}) by substituting the 
canonical pairs $x^1; \dd^1,x^2, \dd^2$ by $x^I; \dd^I,x^{I+1}, \dd^{I+1}$
and the parameters as in the argument.
$$  S^{1234} (2 \frak l^{12},2 \frak l^{34}) =
 S^{23}  (2\underline \ell^1, 2\underline \ell^4) 
 S^{34} ( 2 \underline \ell^2, 2 \underline \ell^4) 
 S^{12} (2\underline \ell^1, 2\underline \ell^3)  S^{23}(2 \underline \ell^2,
2\underline \ell^3) $$
This parameter permutation operator maps equivalent algebras,
$$  S^{1234} (2 \frak l^{12},2 \frak l^{34}) T_4 (u, 2 \frak l^{12},2 \frak l^{34})
= T_4 (u, 2 \frak l^{34},2 \frak l^{12}) S^{1234} (2 \frak l^{12},2 \frak
l^{34}).
$$
$$  S^{1234} (2 \frak l^{12},2 \frak l^{34})\cdot 1 = \pi^{1234} (2 \frak l^{12},2 \frak l^{34})
=
\pi(2\underline \ell^2, 2\underline \ell^3)  \pi(2\underline \ell^1, 2\underline \ell^3)
\pi(2\underline \ell^2, 2\underline \ell^4) \pi(2\underline \ell^1, 2\underline
\ell^4). $$

Now we put 
$$ (2 \underline \ell^3,   2 \underline \ell^4) \to
(2 \underline \ell^1+u , 2 \underline \ell^2+u)  $$
or the analogous with $\ell^1$ and $\ell^2$  permuted. 
With this choice the monodromy $ T_4 $ is composed of two equivalent 
second order factors (\ref{mono2}), in analogy to the discussion about (\ref{pil1}). 

The explicit expression of the permutation coefficient 
$ \pi^{1234} (2 \frak l^{12},2 \frak l^{12} +u)  $ is obtained from (\ref{pi})
by substitutions
$$ \pi(2\underline \ell^I, 2\underline \ell^J +u) = 
\Gamma(2\ell^I_1- 2\ell^J_1-u) \Gamma(2\ell^I_2- 2\ell^J_2-u) 
\frac{\Gamma(2\ell^J_1+1 -2\ell^I_2 +u) }{\Gamma(2\ell^I_1+1 -2\ell^J_2 -u)
} = \pi^{IJ}(u). $$

The reducibility of the representations of the second order evaluation
$T(u; 2\underline \ell^1, 2\underline \ell^2) $ in action on $1$ and
$\psi_0$
is reflected   in the asymptotics of the resulting function in $u \to 0$.

\be \label{pi1234u}
\pi^{1234} (2 \frak l^{12},2 \frak l^{12} +u) = 
\pi^{21}(u) \pi^{11}(u) \pi^{22}(u) \pi^{12}(u).  \ee 
The second and third factors compare to (\ref{pil1}). 
The change in their $u$ asymptotics is related with $2L^1$ or $2L^2 $
being non-negative integer or not.

The first and the last factor compare to (\ref{pil1l2}).   

In the case of generic parameters the asymptotics of $\pi^{1234}$ is 
a constant times $u^{-4} $, 
 the case of $2L^1$ approaching a non-negative integer value
 is distinguished by the asymptotics of the second factor 
 changing from $u^{-2}$ to $-u^{-2}$.

Consider the case of both $2L^1$ and $2L^2$ being non-negative integer
in the product of the first and the fourth factors.
If $2M^{12}$ is generic the asymptotics of their product is 
a constant. If $2M^{12}, 2M^{21}$ are non-negative integer too 
then the asymptotics distinguishes the ordering configurations
and the related type of  representations 
(\ref{L1<MM<L}), (\ref{M12<LL<M}), (\ref{L1=M12}), (\ref{L1=L2=M12})
in the following way.

We write the leading term at $u \to 0$ of the first $\pi^{21}(u)$ and fourth 
$\pi^{12}(u) $ factors
corresponding to the cases 
\begin{enumerate}
\item $ 2L^1 < 2M^{12}, 2M^{12} < 2L^2 $,
$$ \pi^{21} \sim (2L^2-2M^{12}-1)! (-1)^{2L^2-2M^{12}-1} u^{-1}, \ \   
\pi^{12} \sim (2M^{12} -2L^1-1)! (-1^{2M^{12} -2L^1-1} u^{-1}, $$
\item  $2M^{12} < 2L^1, 2L^2 < 2M^{21} $,
$$ \pi^{21} \sim (2L^2-2M^{12}-1)! (2L^1-2M^{12}-1)! (-1)^{2L^1+2L^2+1}, \ \ 
\pi^{12} \sim (-1) u^{-2}, $$
\item $2L^1=2M^{12}, 2L^2=2M^{21}, 2L^2< 2M^{12} $, 
$$ \pi^{21} \sim (-1)^{2M^{12} +2L^2+1} u^{-2}, \ \ \ 
\pi^{12} \sim (2M^{12} -2L^2-1)! (-1 )^{2M^{12} -2L^2-1} u^{-1}, $$
\item $2L^1=2M^{12} =2M^{21} = 2L^2$,
$$ \pi^{21} \sim (-1) u^{-2}, \ \ \ \pi^{12} \sim (-1) u^{-2}. $$
\end{enumerate}

\section {Higher rank permutation coefficients}

\subsection{Beta sequences}

Let us introduce the notation
$$ f(\underline x^{1,i}, \underline x^{2,j}) = det (\underline x^{1,1}, ..., \underline
x^{1,i-1}, \underline x^{1,i+1}, ..., \underline x^{1,n}, \underline x^{2,j}
). $$
The shift operator in $S^1_{k,k+1}$ acts non-trivially on $f(\underline x^{1,i}, \underline x^{2,j})
$ only if the missing vector $\underline x^{1,i}$ is the one with $i=k+1 $.
Only in this case the action contributes to a beta sequence:
\be \label{S1k}
 S^1_{k,k+1}(w) (f(\underline x^{1,k+1}, \underline x^{2,j}))^v =
\int dc c^{w-1} (f(\underline x^{1,k+1}, \underline x^{2,j}) -c f(\underline
x^{1,k}, \underline x^{2,j}) )^v = \ee
 $$
B(w,v+1) (f(\underline x^{1,k}, \underline x^{2,j}))^{-w} (f(\underline
x^{1,k+1}, \underline x^{2,j}))^{u+v}. $$

The shift operator in $S^2_{k,k+1}$ acts non-trivially on $f(\underline x^{1,i}, \underline x^{2,j})
$ only if the  vector $\underline x^{1,j}$ is the one with $j=k $.
Only in this case the action contributes to a beta sequence:
\be \label{S2k}
 S^2_{k,k+1}(w) (f(\underline x^{1,i}, \underline x^{2,k}))^v =
\int dc c^{w-1} (f(\underline x^{1,i}, \underline x^{2,k}) -c f(\underline
x^{1,i}, \underline x^{2,k+1}) )^v = \ee $$
B(w,v+1) (f(\underline x^{1,i}, \underline x^{2,k+1}))^{-w} (f(\underline
x^{1,i}, \underline x^{2,k}))^{u+v}. $$

Consider the parameter array with four of them distinguished on the
positions $(1,i), (1,n), (2,1), (2,i) $ replacing $2\ell^1_{n-i+1},
2\ell^1_1,2\ell^2_n, 2\ell^2_{n-1+1}$ of the standard array, 
$$ 2\ell^1_n, 2\ell^1_{n-1},...,2\bar \ell_b, 2\ell^2_{n-i}, ..., 2\ell^1_2,
2\bar\ell_a; 2\ell_a, 2\ell^2_{n-1}, ..., 2\ell_b, 2\ell^2_{n-i}, ...,
2\ell^2_2, 2\ell^2_1. $$
The permutation $2\bar \ell_a, 2\ell_a$ is represented on the monodromy
by the multiplicative operator
$ S^{12}(2\bar \ell_a-2\ell_a) = (f(\underline x^{1,n},\underline x{^2,1}))^{2\bar \ell_a-2\ell_a}$.
There  are generalizations for the permutations $2\bar \ell_b, 2\ell_a$
and $2\bar \ell_a, 2\ell_b$.

\begin{prop}
A beta sequence of adjacent permutations $ S^{12}_{i,1}$
permuting $2\bar \ell_b, 2\ell_a$ acts as 
\be \label{Si1}
 S^{12}_{i,1} \cdot 1= \pi_{i,1} f(\underline x^{1,i},\underline x^{2,1}))^{2\bar
\ell_b-2\ell_a}. \ee
and a beta sequence of adjacent permutations $ S^{12}_{i,2}$
permuting $2\bar \ell_a, 2\ell_b$ acts as 
\be \label{Si2} S^{12}_{i,2} \cdot 1= \pi_{i,2} f(\underline x^{1,n},\underline x^{2,i}))^{2\bar
\ell_a-2\ell_b}. \ee
The beta sequence permuting $2\bar \ell_b, 2\ell_a$ is constructed as
$$ S^{12}_{i,1} (2\bar \ell_b, 2\bar \ell_a, 2\ell_a) = S^{12}(2\bar \ell_b-2\bar \ell_a)  
\tilde S^{12}_{i,1} (2\bar \ell_b, 2\bar \ell_a, 2\ell_a; 2\frak l\p) S^{12}(2\bar \ell_a-2\ell_a)
$$ $$ = S^{12}(2\bar \ell_b-2\bar \ell_a) 
\stackrel{\rightarrow}S^{1 }_i(2\bar \ell_b, 2\underline \ell^1)
S^1_{i-1,i}(2\bar \ell_b-2\ell_a) \stackrel{\leftarrow}S^{1}_i(2\ell_a,2\underline \ell^1) 
S^{12}(2\bar \ell_a-2\ell_a), $$

$$ \stackrel{\leftarrow}S^{1}_i(2\ell_a,2\underline \ell^1)  = \prod_{n-i-1}^1 S^1_{n-k,
n-k+1}(2\ell^1_{k+1}-2\ell_a), \ \ 
\stackrel{\rightarrow}S^{1 }_i(2\bar \ell_b, 2\underline \ell^1) =
\prod_1^{n-i-1} S^1_{n-k, n-k+1} (2\bar \ell_b-2\ell^1_{k+1} ). $$
The corresponding permutation coefficient depends on the parameters as
$$ \pi^1_i = \stackrel{\rightarrow}\pi^{1}_i B(2\ell_a-2\bar \ell_b,
2\ell^1_{n-i}-2\ell_b+1) \stackrel{\leftarrow}\pi^{1}_i, $$
$$\stackrel{\leftarrow}\pi^{1}_i (2\ell_a,2\underline \ell^1) = \prod_1^{n-i-1}
B(2\ell_a-2\ell^1_{k+1}, 2\ell^1_k-2\ell_{a}+1), $$
$$\stackrel{\rightarrow}\pi^{1}_i (2\bar \ell_b,2\underline \ell^1) = \prod_1^{n-i-1}
B(2\ell^1_{k+1}-2\bar \ell_b, 2\ell^1_k-2\ell^1_{k+1}+1). $$

The beta sequence permuting $2\bar \ell_a, 2\ell_b$ is constructed as
$$ S^{12}_{i,2} (2\bar \ell_a, 2\ell_a, 2\ell_b) = S^{12}(2 \ell_b-2\ell_a)
\tilde S^{12}_{i,2} (2\bar \ell_a, 2\ell_a, 2\ell_b) S^{12}(2 \bar \ell_a-2\ell_a)
=$$ $$ S^{12}(2 \ell_b-2\ell_a)
\stackrel{\leftarrow}S^{2}_i(2\ell_b; 2\underline \ell^2) S^2_{i-1,i}(2\bar \ell_a-2\ell_b)
\stackrel{\rightarrow}S^{2}_i(2\bar \ell_a; 2\underline \ell^2) S^{12}(2 \bar
\ell_a-2\ell_a), $$

$$
\stackrel{\rightarrow}S^{2}_i(2\bar \ell_a; 2\underline \ell^2) =
\prod_{k=i-2}^1S^2_{k,k+1}(2\bar \ell_a- 2\ell^2_{n-k}), \ \ 
\stackrel{\leftarrow}S^{2}_i(2\ell_b; 2\underline \ell^2) = \prod_{k=1}^{i-2}
S^2_{k,k+1}(2\ell^2_{n-k}-2\ell_b ).  $$
The corresponding permutation coefficient depends on the parameters as
$$\pi_{i,2}(2\bar \ell_b, 2\ell_a) = \stackrel{\leftarrow}\pi_{i,2} (2\ell_b; 2\underline \ell^2) B(2\ell_b-2\bar
\ell_a, 2\bar \ell_a -2\ell_{n-i+2}+1) \stackrel{\rightarrow}\pi_{i,2} (2\bar \ell_a; 2\underline
\ell^2),
$$
$$ \stackrel{\leftarrow}\pi_{i,2} (2\ell_b; 2\underline \ell^2) = \prod_{k=1}^{i-2}
B(2\ell_b-2\ell^2_{n-k}, 2\ell^2_{n-k}- 2\ell^2_{n-k+1} +1), $$

$$ \stackrel{\rightarrow}\pi_{i,2} (2\bar \ell_a; 2\underline \ell^2) = 
\prod_{k=1}^{i-2} B(2\ell^2_{n-k}-2\bar \ell_a, 2\bar \ell_a -
2\ell^2_{n-k+1 }+1). $$

\end{prop}

{\it Proof:}

We start with
$$ S^{12}(2\bar \ell_a-2\ell_a) \cdot 1 = (f(\underline x^{1,n}, \underline
x^{2,1}) )^{2\bar \ell_a-2\ell_a} $$
and use (\ref{S1k}) to calculate  the action of the ordered product 
$\stackrel{\leftarrow}S^{1}_i $ on this function to obtain
$$ \stackrel{\leftarrow}\pi^{1}_i (2\ell_a,2\underline \ell^1)  
(f(\underline x^{1,n}, \underline x^{2,1}) )^{2\bar \ell_a -2\ell^1_2}
\prod_{k=2}^{n-i-1} f(\underline x^{1, i+k}, \underline
x^{2,1}) )^{2\ell^1_{n-i+1-k} -2\ell^1_{n-i+2-k} }
$$ $$ 
(f(\underline x^{1,i+1}, \underline x^{2,1}) )^{2 \ell^1_{n-i} - 2\bar
\ell_b}. 
$$ 
After this action the permuted parameter array is
$$ 2\ell^1_n, 2\ell^1_{n-1},...,2\bar \ell_b,  2\ell_a, 2\ell^2_{n-i}, ...,
2\ell^1_2; 
 2\bar \ell_a, 2\ell^2_{n-1}, ..., 2\ell_b, 2\ell^2_{n-i}, ...,
2\ell^2_2, 2\ell^2_1. $$

The action of the next adjacent permutation operator $ S^1_{i-1,i}(2\bar \ell_b-2\ell_a)  $ permuting $2\bar
\ell_b,2\ell_a $ results in the multiplication of the
latter by 
$$ B(2\ell_b-2\bar \ell_a, 2\bar \ell_a -2\ell_{n-i+2}+1)  
(f(\underline x^{1,i}, \underline x^{2,1}) )^{2\bar \ell_b -2\ell_a}.   
$$ 
The action of $ \stackrel{\rightarrow}S^{1 }_i(2\bar \ell_b, 2\underline \ell^1)  $ on the
result leaves the last factor unchanged. The factors involving 
$\underline x^{1,i+k}, k=1,..., n-i-1 $ are canceled. 
After this action the permuted parameter array is
$$ 2\ell^1_n, 2\ell^1_{n-1},...,2 \ell_a, 2\ell^1_{n-i}, ...,
2\bar \ell_b; 
 2\bar \ell_a, 2\ell^2_{n-1}, ..., 2\ell_b, 2\ell^2_{n-i}, ...,
2\ell^2_2, 2\ell^2_1. $$
We obtain 
$$ \stackrel{\rightarrow}S^{1 }_i(2\bar \ell_b, 2\underline \ell^1)
S^1_{i-1,i}(2\bar \ell_b-2\ell_a) \stackrel{\leftarrow}S^{1}_i(2\ell_a,2\underline \ell^1) 
S^{12}(2\bar \ell_a-2\ell_a) \cdot 1 = $$ $$
\stackrel{\rightarrow}\pi_{i,2} (2\bar \ell_a; 2\underline \ell^2)
B(2\ell_b-2\bar \ell_a, 2\bar \ell_a -2\ell_{n-i+2}+1) 
\stackrel{\leftarrow}\pi^{1}_i (2\ell_a,2\underline \ell^1) $$ $$
(f(\underline x^{1,i}, \underline x^{2,1}) )^{2\bar \ell_b -2\ell_a}
(f(\underline x^{1,n}, \underline x^{2,1}) )^{2\bar \ell_a -2\bar \ell_b}.
$$
Now the action of $S^{12}(\bar \ell_b- 2\bar \ell_a)$ cancels the last
factor. We obtain the asserted expression and the wanted 
 parameter array 
$$ 2\ell^1_n, 2\ell^1_{n-1},...,2 \ell_a, 2\ell^2_{n-i}, ...,
2\bar \ell_a; 
 2\bar \ell_b, 2\ell^2_{n-1}, ..., 2\ell_b, 2\ell^2_{n-i}, ...,
2\ell^2_2, 2\ell^2_1. $$

The action  $S^{12}_{i,2} (2\bar \ell_b, 2\ell_a) \cdot 1$
is calculated in analogy.

\qed

$\tilde S^{12}_{i,1}$ and the corresponding permutation factor depend besides of the parameters indicated in the
argument also on the parameters $2\ell^1_{i-1}, ..., 2\ell^1_2$.
$\tilde S^{12}_{i,2}$ depends besides of the parameters indicated in the
argument also of the parameters $2\ell^2_{n-1}, ..., 2\ell^2_{i+1}$.

The expressions of the permutation operators $ S^{12}_{i,1}$ and $
S^{12}_{i,2}$ are related by replacing the superscript $1 \to 2$, reflecting
the arrows, substituting of the running index $k +1 \to n-k$, the
distinguished parameters  
$2\ell_a \to 2 \bar \ell_a, 2 \bar \ell_b \to 2\ell_b $ and modifying the
adjacent permutation operator factors 
$ S^2_m (2 L) $ additionally by reversing the sign in the argument and
diminishing the subscript by 1.  The expressions of the permutation
coefficients $ \pi^{12}_{i,1}$ and $\pi^{12}_{i,2}$ are related by te same
substitutions and by modifying the arguments of the factors 
$B(2L, 2M+1)$ additionally by reversing the signs of $2L$ and $2M$.

\subsection{1 to 1 beta sequences}

We intend to construct the 1 to 1 beta sequences of adjacent permutations
representing the permutation $\sigma^{12}_i$ exchanging
$2\ell^1_i \leftrightarrow 2\ell^2_i $.
\be \label{Si} S^{12}_i T(u,2\frak l) = T(u,2\sigma^{12}_i \frak l) S^{12}_i,
\ \ \ 
 S^{12}_i \cdot 1 = \pi^{12}_i \cdot 1 .\ee

Both monodromies have $1$ as a highest weight vector. The weight functions
coincide. In the case of generic parameters the representations are of
highest weight and they are equivalent. The permutation coefficients
 $ \pi^{12}_i$ carry the information about the parameter values of the cases
of non-equivalence.  

We notice that, disregarding sign,
\be \label{fxj} f(\underline x^{1,j}, \underline x^{2,j }) = 1. \ee

The last proposition implies the result for the cases $ i=1$ and $i=n$
immediately.

$$ S^{12}_1 = S^{12}( 2 \ell_a-2 \ell_b) \tilde S^{12}_{1,2} (2\ell^1_1,
2\ell^2_n, 2\ell^2_1) S^{12} (2\bar \ell_a-2\ell_a), $$
\be \label{S1} \pi^{12}_1 = \pi^{12}_{1,2} (2\ell^1_1, 2\ell^2_n,2\ell^2_1),
\ee
$$ S^{12}_n = S^{12} (2\bar \ell_a-2\bar \ell_b)
\tilde S^{12}_{n,1} (2\ell^1_n,2\ell^1_1; 2\ell^2_n) S^{12} (2\bar
\ell_a-2\ell_a),$$
\be \label{Sn} \pi^{12}_n = \pi^{12}_{n,1} (2\ell^1_n,2\ell^1_1, 2\ell^2_n).
\ee

To formulate the operators obeying (\ref{Si}) for the remaining values of $i$ we
rewrite (\ref{Si1}) and (\ref{Si2}). 
\be \label{tildeSi1} 
\tilde S^{12}_{i,1} (2\bar \ell_b,2\bar \ell_a, 2\ell_a) (f(\underline x^{1,n}, \underline x^{2,1}) )^{2\bar \ell_a -2
\ell_a} = $$ $$ \pi^{12}_{i,1} 
(f(\underline x^{1,n}, \underline x^{2,1}) )^{2 \bar \ell_b -2
\bar \ell_b}
(f(\underline x^{1,i}, \underline x^{2,1}) )^{2\bar \ell_b -2 \ell_a},
\ee
\be \label{tildeSi2}
\tilde S^{12}_{i,2} (2\bar \ell_a, 2\ell_a, 2\ell_b) (f(\underline x^{1,n}, \underline x^{2,1}) )^{2\bar \ell_a -2
\ell_a} = \pi^{12}_{i,2}
(f(\underline x^{1,n}, \underline x^{2,1}) )^{2 \ell_a -2 \ell_b}
(f(\underline x^{1,n}, \underline x^{2,i}) )^{2\bar \ell_a -2 \ell_b}.
\ee

\begin{prop}
 1 to 1 beta sequences of permutation operators $ S^{12}_i$ acting as in
(\ref{Si}) are for $i=2, ..., n-1$

$$ S^{12}_i (2 \frak l) = S^{12}(2\ell^2_n-2\ell^1_1) 
\tilde S^{12}_{n-i+1, 1}(2\ell^2_n 2\ell^1_1,2\tilde \ell_i, 2\bar
\ell_i) $$ $$
S^{12}(2\ell^1_1-2\tilde \ell_i) \tilde S^{12}_{n-i+1,2}(2\ell^2_n,2\bar \ell_i,
2\ell^2_n, 2\tilde \ell_i)  S^{12}(2\ell^1_i-2\ell^1_1) \tilde S^{12}_{n-i+1, 1}(2\bar \ell_i, 2\ell^1_1,
2\ell^2_n, 2\tilde \ell_i ) S^{12}(2\ell^1_1-2\ell^1_n), $$
$$
\pi^{12}_i (2 \frak l) =
\pi^{12}_{n-i+1, 1}(2\ell^2_n 2\ell^1_1,2\tilde \ell_i, 2\bar \ell_i)
\pi^{12}_{n-i+1,2} ((2\ell^2_n,2\bar \ell_i,2\ell^2_n, 2\tilde \ell_i)
\pi^{12}_{n-i+1, 1}(2\bar \ell_i, 2\ell^1_1,
2\ell^2_n, 2\tilde \ell_i ). 
$$

The cases $i=1, n $ are given above (\ref{S1}), (\ref{Sn}).

\end{prop}

{\it Proof:}

We start with the operation permuting  $2\ell^1_i, 2\ell^2_n$ (\ref{Si1}),
$$ S^{12}_{n-i+1, 1}(2\bar \ell_i, 2\ell^1_1,2\ell^2_n, 2\tilde \ell_i )
\cdot 1 = \pi^{12}_{n-i+1, 1}(2\bar \ell_i, 2\ell^1_1,
2\ell^2_n, 2\tilde \ell_i ) \ 
(f(\underline x^{1,n-i+1}, \underline x^{2,1}))^{\bar\ell_i-2\ell^2_n}.  
$$
This operation changes the parameter array to 
$$ 2\ell^1_n, 2\ell^1_{n-1},...,2\ell^2_n,  2\ell^1_{i-1},  ...,
2\ell^1_2, 2 \ell^1_1; 2\ell^1_i, 2\ell^2_{n-1},..., 2\ell^2_i,
2\ell^2_{i-1}, ...,2\ell^2_2, 2\ell^2_1. $$
The action of $\tilde S^{12}_{n-i+1,2}$ on
$(f(\underline x^{1,n-i+1}, \underline x^{2,1}))^{\ell^1_i-2\ell^2_n}$ 
is obtained as the result of (\ref{tildeSi2}) however for the parameter
configuration
$$ ..., 2\ell^1_i; 2\ell^2_n, 2\ell^2_{n-1},..., 2\ell^1_i,
2\ell^2_{i-1},...$$
with the result
$$ \tilde S^{12}_{n-i+1,2}(2\ell^1_n,2\ell^1_i,2\ell^2_n, 2\ell^2_i)  
(f(\underline x^{1,n-i+1}, \underline x^{2,1}))^{\ell^1_i-2\ell^2_n} =
$$ $$
\pi^{12}_{n-i+1,2}(2\ell^1_n,2\ell^1_i,2\ell^2_n, 2\ell^2_i) 
(f(\underline x^{1,n-i+1}, \underline x^{2,1}))^{2\ell^2_i-2\ell^2_n}
(f(\underline x^{1,n-i+1}, \underline x^{2,n-i+1}))^{\ell^1_i-2\ell^2_i}.
$$ 
According to (\ref{fxj}) the last factor can be replaced by 1. 
This operation changes the parameter array to 
$$ 2\ell^1_n, 2\ell^1_{n-1},...,2\ell^2_n,  2\ell^1_{i-1},  ...,
2\ell^1_2, 2 \ell^1_1; 2\ell^2_i, 2\ell^2_{n-1},..., 2\ell^1_i,
2\ell^2_{i-1}, ...,2\ell^2_2, 2\ell^2_1. $$
The following $ S^{12}(2\ell^1_1-2\ell^2_i) $ permutes 
$  2 \ell^1_1, 2\ell^2_i$. This operation does not add a parameter dependent
factor and results in the multiplication by 
$$  (f(\underline x^{1,n}, \underline x^{2,1}))^{2\ell^1_1-2\ell^2_i}. $$
 The  action of $\tilde S^{12}_{n-i+1, 1}(2\ell^2_n 2\ell^1_1,2\tilde \ell_i, 2\bar
\ell_i) $ is calculated according to (\ref{tildeSi1}), it acts non-trivially on the
latter factor only. The result is  multiplied by
$$
\pi^{12}_{n-i+1,1}(2\ell^2_n 2\ell^1_1,2\tilde \ell_i, 2\bar
\ell_i) 
(f(\underline x^{1,n}, \underline x^{2,1}))^{2\ell^1_1-2\ell^2_n} 
(f(\underline x^{1,n-i+1}, \underline x^{2,1}))^{2\ell^2_n-2\ell^2_i}.
$$

The factors with $\underline x^{1,n-i+1}$ in the argument cancel. 
This operation changes the parameter array to 
$$ 2\ell^1_n, 2\ell^1_{n-1},...,2\ell^2_i,  2\ell^1_{i-1},  ...,
2\ell^1_2, 2 \ell^2_n; 2\ell^1_1, 2\ell^2_{n-1},..., 2\ell^1_i,
2\ell^2_{i-1}, ...,2\ell^2_2, 2\ell^2_1. $$
The last operation $S^{12} (2\ell^2_n-2\ell^1_1)$
permutes $ 2 \ell^2_n, 2\ell^1_1$ leading to the wanted parameter array.
Its action cancels the factors  with $\underline x^{1,n-i+1}$ in the argument
without adding a parameter dependent factor. 

\qed 

The results provide the basis of investigating the types of
representations with a highest weight of the order $N$ evaluations
of the Yangian algebra of the type $g\ell(n)$ by considerations in analogy to
sect. 8.4. The main features appear in the second order evaluation. They can
be analysed by the fourth order evaluation $T_4(v; 2\underline \ell^1, 
2\underline \ell^2, 2\underline \ell^3, 2\underline \ell^4) $ and the $u$
dependence of the permutation coefficient for 
$$  2\underline \ell^1, 2\underline \ell^2, 2\underline \ell^3, 2\underline \ell^4
\to 2\underline \ell^3, 2\underline \ell^4, 2\underline \ell^1, 2\underline \ell^2
$$
with $(2\underline \ell^3, 2\underline \ell^4) = (2\underline \ell^1+u, 2\underline
\ell^2+u ) $.

\vspace{1cm}

\noindent
{\bf Acknowledgments.}

\noindent
The author is grateful to A. Molev for useful discussions. 

\vspace{1cm}


\end{document}